\documentclass[12pt,a4paper,reqno]{amsart}  

\usepackage[headings]{fullpage} 

\usepackage{graphicx}
\usepackage{float}  
\usepackage{enumerate}   
\usepackage{datetime} 
\usepackage{url}  
\usepackage[english]{babel}

\usepackage{bm} 
\usepackage{mathrsfs}
  
\usepackage[full]{textcomp} 
\usepackage{newtxtext} 
\usepackage{cabin} 
\usepackage{zlmtt}
\usepackage[bigdelims,vvarbb]{newtxmath} 
\usepackage[cal=boondoxo]{mathalfa} 

\usepackage{microtype} %
\newcommand{\eps}{\varepsilon}
\newcommand{\dx}{\mathrm{d}} 
\newcommand{\A}{\mathcal{A}} 
\newcommand{\B}{\mathcal{B}} 
\newcommand{\CC}{\mathcal{C}}

\newcommand{\C}{\mathbb{C}}
\newcommand{\Q}{\mathbb{Q}}
\newcommand{\N}{\mathbb{N}}
\newcommand{\R}{\mathbb{R}} 

\newcommand{\Z}{\mathbb{Z}}
\newcommand{\G}{\mathfrak{G}}

\newcommand{\LL}{\mathscr L}

\newcommand{\Odi}[1]{\Odip{}{#1}}

\newcommand{\Odip}[2]{\mathcal{O}_{#1}\left(#2\right)}

\newcommand{\odi}[1]{{o}\left(#1\right)}

\allowdisplaybreaks

\renewcommand{\qedsymbol}{$\square$}
\newenvironment{Proof}[1][Proof]{\par\noindent\textbf{#1.}~}
{\hfill\qedsymbol\smallskip\par}

\newtheoremstyle{sltheorems}
{10pt}
{6pt}
{\slshape}
{}
{\bfseries}
{.}
{.5em}
{\thmname{#1}\thmnumber{ #2}\thmnote{ (#3)}}
\theoremstyle{sltheorems} 
\newtheorem{Theorem}{Theorem}
\newtheorem{Corollary}{Corollary} 
\newtheorem{Lemma}{Lemma} 
\newtheorem{Proposition}{Proposition} 

\newtheoremstyle{remark}
{10pt}
{6pt}
{\rm} 
{}
{\bfseries}
{.}
{.5em}
{\thmname{#1}\thmnumber{ #2}\thmnote{ (#3)}}
\theoremstyle{remark} 
\newtheorem{Remark}{Remark}

\usepackage{caption} 
\allowdisplaybreaks
\makeatletter
\def\subsubsection{\@startsection{subsubsection}{3}%
\z@{.3\linespacing\@plus.5\linespacing}{-.5em}%
 {\normalfont\bfseries}}
\makeatother

\newcommand{\bound} {10^7}

\begin{document}
\let \oldsection \section
\def \newsection {\vspace{14pt plus 3pt}\oldsection}
\renewcommand{\section}{\newsection}

\title[A fast algorithm to compute Ramanujan-Deninger's gamma-function]{A fast
algorithm to compute the Ramanujan-Deninger \\
gamma-function and some number-theoretic applications}
\author{Alessandro Languasco and Luca Righi}

\begin{abstract}  
We introduce a fast algorithm to compute the Ramanujan-Deninger gamma function
and its logarithmic derivative at positive  values.
Such an algorithm allows us to greatly extend the numerical investigations about
the Euler-Kronecker constants $\G_q$, $\G_q^+$ and
$M_q=\max_{\chi\ne \chi_0} \vert L^\prime/L(1,\chi)\vert$,
where $q$ is an odd prime, $\chi$ runs over the primitive Dirichlet characters 
$\bmod\ q$, $\chi_0$ is the principal Dirichlet character $\bmod\ q$ and 
$L(s,\chi)$ is the Dirichlet $L$-function associated to $\chi$.
Using such algorithms we  obtained that
$\G_{50 040  955 631} =-0.16595399\dotsc$ and
$\G_{50 040  955 631}^+ =13.89764738\dotsc$
thus getting a new negative value for $\G_q$.

Moreover we also computed $\G_q$, $\G_q^+$ and $M_q$ 
for every  prime $q$, $10^6< q\le \bound$, thus extending the results in 
\cite{Languasco2021a}.
As a consequence we obtain that both $\G_q$ and $\G_q^+$ 
are positive for every odd prime $q$ up to $\bound$
and  that
\(
\frac{17}{20} \log \log q< M_q < \frac{5}{4} \log \log q
\)
for every  prime $1531 < q\le \bound$. In fact the lower bound 
holds true for  $q>13$. 
The programs used and the results here described are collected 
at the following address
\url{http://www.math.unipd.it/~languasc/Scomp-appl.html}.
\end{abstract} 
\maketitle 
\mbox{}\vskip-1.75truecm
\section{Introduction} 
We introduce a fast algorithm to compute the 
Ramanujan-Deninger gamma function and its
logarithmic derivative at positive values.
We then use such a new algorithm to efficiently compute
 $L^\prime/L(1,\chi)$, where  $\chi$ runs over the non-principal
 primitive Dirichlet characters $\bmod\ q$,
$q$ is an odd prime and $L(s,\chi)$ is the Dirichlet 
$L$-function associated to $\chi$.
Such a quantity is involved in several interesting
number-theoretic problems like the evaluation of
the Euler-Kronecker constants $\G_q$ for the cyclotomic
field $\Q(\zeta_q)$, $\zeta_q$ being a $q$-root of unity,
the analogous problem for $\G_q^+$ attached to 
$\Q(\zeta_q+\zeta_q^{-1})$, the maximal
real subfield of $\Q(\zeta_q)$,
and the study of the extremal values of $M_q=\max_{\chi\ne \chi_0} \vert L^\prime/L(1,\chi)\vert$. 
We will give a  detailed description of such problems in Section \ref{sect-applications}.

Following Deninger's notation in \cite{Deninger1984}, we introduce now the 
functions we will work on. The main object is the \emph{Ramanujan-Deninger} Gamma function
$\Gamma_1(x):=\exp(R(x))$, $x>0$, where 
\begin{equation*}
R(x) := - \frac{\partial^2}{\partial s^2} \zeta(s,x) \vert _{s=0}\, , 
\end{equation*}
 $\zeta(s,x)$ is the Hurwitz zeta function,
$\zeta(s,x) = \sum_{n=0}^{+\infty}(n+x)^{-s}$ for $\Re(s)>1$
and it is  meromorphically extended to
$s\in \C\setminus\{1\}$.
We recall that $\zeta(s,1)$ is Riemann's zeta-function $\zeta(s)$.
Using eq.~(2.3.2) of \cite{Deninger1984}, 
the $R$-function can be expressed for every $x>0$ by
\begin{align}
\label{R-equiv}
R(x) &= -\zeta^{\prime\prime}(0) - S(x),\\
\label{S-def}
S(x) 
&:=
2 \gamma_1 x +(\log x)^{2} 
+ 
\sum_{k=1}^{+\infty}
\Bigl(
\bigl(\log (k+x)\bigr)^{2} - (\log k)^{2} -2x \frac{\log k}{k}
\Bigr),
\end{align} 
where
\begin{equation} 
\label{gamma1-zetasecond-def}
\gamma_{1} =
\lim_{n\to+\infty}
\Bigl(
\sum_{k=1}^{n}
\frac{\log k}{k} - \frac{(\log n)^2}{2}
\Bigr) ,
\quad
\zeta^{\prime\prime}(0)= \frac12 \Bigl(-(\log 2\pi)^2-\frac{\pi^2}{12} +
\gamma_1 + \gamma^2 \Bigr)
\end{equation}
and $\gamma$ is the Euler-Mascheroni constant.
We introduced the $S$-function because in the applications we will see
in Section \ref{sect-applications} below the constant term
$\zeta^{\prime\prime}(0)$ will play no role
and hence we may focus our attention
just on the $S$-function. 
We also have $S(1)=0$ and $R(1) = -\zeta^{\prime\prime}(0)$.

A key point to be able to obtain the following results is that $R(x)$ is the
unique solution
in $(0,+\infty)$ of the difference equation $R(x+1) = R(x) + (\log x)^2$, 
with initial condition $R(1) = -\zeta^{\prime\prime}(0)$, which is convex in
some interval $(A,+\infty)$, $A>0$, see Theorem 2.3 of Deninger
\cite{Deninger1984}.
As a consequence the $S$-function too verifies a difference equation:
\begin{equation}
\label{difference-S}
S(x+1) = S(x) - (\log x)^2 \quad \textrm{for every} \ x>0.
\end{equation}
Another important ingredient to be used to derive our results is the
following alternative definition of $S(x)$, $x>0$, which 
is implicitly contained in eq.~(2.12) of Deninger \cite{Deninger1984}.
It is the analogue for $S$ of \emph{Plana's integral} for $\log \Gamma$,
where $\Gamma(s)$, $s\in \C\setminus(-\N)$, is Euler's Gamma function:
\begin{equation}
\label{S-alt-def} 
S(x) =
2 \int_{0}^{+\infty} 
\Bigl( (x-1) e^{-t} + \frac{e^{-xt} - e^{-t}}{1-e^{-t}} \Bigr)
\frac{\gamma + \log t}{t} \ \dx t.
\end{equation}
We introduce now the first derivative of $R(x)$, namely
\begin{equation}
\label{psi1-def}
\psi_{1}(x):=\frac12 R^\prime(x)  =
\frac12 \frac{\Gamma_1^\prime}{\Gamma_1} (x).
\end{equation}
The factor $1/2$ in \eqref{psi1-def} is needed because Deninger, in its
definition of $R(x)$,
used an extra factor $2$ and we need now to remove it to connect $\psi_1$
with the results proved by other authors.
We also recall that generalised $\psi$-functions of this kind occur in
Ramanujan's second notebook, see \cite[Chapter 8, Entry 22]{Berndt1985}.
Differentiating \eqref{S-def}
we have
\begin{equation*}
\psi_1(x) =
  -\gamma_1 - \frac{\log x}{x}
 - \sum_{k=1}^{+\infty}
\Bigl(
\frac{\log (k+x)}{k+x} - \frac{\log k}{k}
\Bigr) 
\end{equation*}
for  $x>0$. We also define
\begin{equation}
\label{T-def}
T(x) 
:=
\gamma_1 + \psi_{1}(x)
\end{equation} 
 so that  $T(1) =0$.
 As before, we introduced the function $T(x)$ because the constant term
$\gamma_1$ in the definition of $\psi_1$ will play no
 role in the applications contained in Section \ref{sect-applications} below.
Moreover, since  
$S^\prime(x) = 2(\gamma_1 - T(x))$, 
it follows that
\begin{equation}
\label{TS-link-difference}
T(x)
= \gamma_1 - \frac12 S^\prime(x) 
\quad
\textrm{and}
\quad
 T(x+1) = T(x) + \frac{\log x}{x}
\quad \textrm{for every} \ x>0.
\end{equation}
In our applications, see Section \ref{sect-applications} below, we will need to
evaluate $S$ or $T$ at some rational points
contained in $(0,1)$. A  possible solution, used in \cite{Languasco2021a}, 
is to use the \texttt{intnum} and \texttt{sumnum}
functions of PARI/GP \cite{PARI2021} to numerically evaluate \eqref{S-def},
\eqref{S-alt-def} and \eqref{T-def}. 
Here we show how to largely reduce the cost of this computation by introducing
a new algorithm to obtain such quantities.
Denoting as $\lceil y \rceil$ 
the least integer greater than or equal to 
$y\in \R$ and 
defining the $m$-th \emph{harmonic} number as
\begin{equation}
\label{harmonic-def}
 H_m =\sum_{j=1}^{m}\frac{1}{j},
\end{equation}
where $m\in\N$, $m\ge1$,
our starting point is the following  Theorem  \ref{main-thm-S}
in which   formula \eqref{S-new-series} was  first proved\footnote{Pay attention to the 
fact that the Deninger $S(x)$-function defined in
\eqref{R-equiv}-\eqref{S-def} is equal to 
$ -2 \log(\Gamma_1(x))$ as defined 
in Proposition 1 of Dilcher \cite{Dilcher1994}.} in Dilcher
\cite[eq.~(2.14)]{Dilcher1994}. 
\begin{Theorem}
\label{main-thm-S}
Let $x\in (0,2)$. Then
\begin{equation}
\label{S-new-series}
S(x) = 
-
2 \gamma_1 (1-x) 
+
2 \sum_{k=2}^{+\infty} \frac{\zeta(k) H_{k-1}+\zeta^\prime(k)}{k} (1-x)^{k}, 
\ \end{equation}
where 
$\gamma_1$ is defined as in \eqref{gamma1-zetasecond-def},
$H_{k}$ is defined as in \eqref{harmonic-def},
$\zeta(\cdot)$ is the Riemann zeta-function and $\zeta^\prime(\cdot)$ is its first
derivative.
Moreover, letting $x\in (0,1)\cup(1,2)$, $n\in \N$, $n\ge 1$ be fixed, 
and $r_S (x,n)\in \N$,
\begin{equation}
\label{r-estim-S}
r_S (x,n) :=
\Bigl\lceil
\frac{(n+2) \log 2 + \vert\log (1-\vert 1-x\vert) \vert
}{\vert \log\vert 1-x \vert \vert}
\Bigr\rceil
-1,
\end{equation}
 we have that there exists
$\theta=\theta(x)\in (-1/2,1/2)$ such that
\begin{equation}
\label{S-new-trunc-series}
S(x) = 
-
2 \gamma_1 (1-x) 
+
2 \sum_{k=2}^{r_S (x,n)}\frac{\zeta(k) H_{k-1}+\zeta^\prime(k)}{k} (1-x)^{k}
+
\
\vert \theta\vert 2^{-n}. 
\end{equation}
\end{Theorem}

We immediately remark that \eqref{S-new-series} is the Taylor series centred at $1$
of $S(x)$; in particular this implies  that
\[
S^\prime(1) = 2\gamma_1
\quad
\text{and}
\quad
S^{(k)} (1)=
2(-1)^k (k-1)! \bigl( \zeta(k) H_{k-1} + \zeta^\prime(k)\bigl) 
\quad (k\in \N, k\ge 2).
\] 

Using Theorem \ref{main-thm-S} and \eqref{R-equiv}, for $x\in(0,2)$ we trivially have 
\begin{align*} 
R(x) &=
 -\zeta^{\prime\prime}(0) 
+
2(1-x)\gamma_1 
-
2 \sum_{k=2}^{+\infty} \frac{\zeta(k) H_{k-1}+\zeta^\prime(k)}{k} (1-x)^{k}
\end{align*} 
and its corresponding truncated version.
We already remarked that \eqref{S-new-series} is equation (2.14) of Dilcher
\cite{Dilcher1994} but in Theorem \ref{main-thm-S} we will prove it in a
different way, \emph{i.e.}, starting from \eqref{S-alt-def}, which in fact 
reveals that such an argument can be used for any function having an integral 
representation of Plana's type like the one in \eqref{S-alt-def}.
Formulae \eqref{r-estim-S}-\eqref{S-new-trunc-series} of Theorem \ref{main-thm-S} 
are new.

Recalling \eqref{difference-S}, the fact that Theorem \ref{main-thm-S} holds
for every $x\in(0,2)$ means that every value of $S(x)$, $x\in(0,1)$, can be
computed in two different ways.\footnote{We remark
that the size of the convergence interval of the series in the right hand side
of \eqref{S-new-series} can be doubled by isolating
the Taylor series at $1$ of $(\log x)^2$ and using the estimates
on $\vert \zeta(n)-1 \vert$ of Lemma \ref{elementary-estim} below. We do not
insert such an idea here,
since the computation of such an extra-factor $(\log x)^2$ leads, in our
practical application,
to a longer total running time.} Moreover it is clear 
that $r_S (x,n)$ becomes larger as $\vert 1-x\vert$ increases.  
Hence, remarking that in our applications of Section \ref{sect-applications} 
we are mainly interested in $x\in(0,1)$,   
if $x\in (1/2,1)$ we will directly compute $S(x)$ using Theorem \ref{main-thm-S}
while for $x\in (0,1/2)$ we will shift the problem using \eqref{difference-S}
and use Theorem \ref{main-thm-S} in $(1,3/2)$. 
In the following we will refer to this procedure as the \emph{shifting trick}.
Such an argument leads to the following two corollaries.
\begin{Corollary}
\label{corollary-S}
Let $x\in (0,1/2)$. We have that 
\begin{equation}
\label{S-periodic-series}
S(x) =
(\log x)^2 
+
2 \gamma_1 x
+
2 \sum_{k=2}^{+\infty}(-1)^{k}\frac{\zeta(k) H_{k-1} +\zeta^\prime(k)}{k} x^{k}.
\end{equation}
 Letting further $n\in \N$, $n\ge 1$ be fixed
and  $r^\prime_S (x,n) \in \N$, 
\begin{equation*}
r^\prime_S(x,n):=r_S (1+x,n) = 
\Bigl\lceil
\frac{(n+2) \log 2 +\vert \log (1-x) \vert }{
\vert \log x \vert}
\Bigr\rceil
-1,
\end{equation*}
where $r_S (u,n)$ is defined in Theorem \ref{main-thm-S}, %
we have that there
exists $\eta=\eta(x)\in (-1/2,1/2)$ such that
\begin{equation}
\label{S-periodic}
S(x) =
(\log x)^2 
+
2 \gamma_1 x
+
2 \sum_{k=2}^{r^\prime_S (x,n)}(-1)^{k}\frac{\zeta(k) H_{k-1} +\zeta^\prime(k)}{k} x^{k}
+
\
\vert \eta\vert 2^{-n}.
\end{equation}
\end{Corollary}
Recalling Remark (2.6) of Deninger \cite[page 176]{Deninger1984}, 
\eqref{R-equiv} and \eqref{gamma1-zetasecond-def} we also obtain
\begin{equation}
\label{S-value-12}
S\Bigl(\frac12\Bigr)
=
- R\Bigl(\frac12\Bigr)- \zeta^{\prime\prime}(0) 
=
\frac12 (\log \pi)^2
+
\frac{\pi^2}{24}
- \frac{\gamma_1+\gamma^2}{2}.
\end{equation}
\begin{Remark}
As a matter of curiosity, since we are aware of the fact that much faster
algorithms exist to compute $\gamma_1$, see Johansson-Blagouchine
\cite{JohanssonB2019}, we remark that
evaluating twice $S(1/2)$ using Theorem \ref{main-thm-S}
(the first time directly and the second as $S(3/2)+(\log 2)^ 2$), we have, by
subtracting
such formulae, that the summands having even indices vanish; thus we obtain
\[
\gamma_1
= 
- \frac12
(\log 2)^2 
+
\sum_{\ell=1}^{n/2+1} \frac{ \zeta(2\ell+1) H_{2\ell} + \zeta^\prime(2\ell+1)
}{(2\ell+1)4^{\ell}} 
+
\tau 2^{-n+1},
\]
for some $\tau\in(-1/2,1/2)$.
Such a result,
similar to equation (3.10) of Dilcher \cite{Dilcher1994},
 allows us to fast compute $\gamma_1$ with a precision of $n$ bits using about
$n/2$ summands. For example, using PARI/GP, we got $\gamma_1$ with  a precision
of $1\,000$ decimal digits within 
$1$ minute and $1$ seconds of computation time on a Dell
 OptiPlex-3050 machine (equipped with an Intel i5-7500 processor, 3.40GHz, 
16 GB of RAM and running Ubuntu 18.04.5).
\end{Remark}

Combining \eqref{S-new-trunc-series}, \eqref{S-periodic} and \eqref{S-value-12},
we obtain a very fast way of computing $S(x)$ for every $x>0$. We will see
more about this in Section \ref{sect-implementation} but we also summarise 
the situation in the following
\begin{Corollary}
\label{S-computation}
We use the notation introduced in Theorem \ref{main-thm-S} and Corollary
\ref{corollary-S}. 
Moreover, for every $x>0$, we denote
as $\lfloor x \rfloor$ the integral part of $x$ and as $\{x\} = x - \lfloor x
\rfloor$ the fractional part of $x$. Hence we obtain:
\begin{enumerate}[i)]
\item
 $S(1)=S(2)=0$ and  $S(m) = -\sum_{k=2}^{m-1} (\log k)^2$
for every $m\in \N$, $m\ge 3$;
\item
for $x>1$, $x\not \in \N$, we compute $S(x)$ as $S(x) = S(\{x\}) -
\sum_{k=0}^{\lfloor x \rfloor - 1} (\log (\{x\}+k))^2$;
\item
$S(1/2 ) =  (\log \pi)^2/2
+
\pi^2/24
- (\gamma_1+\gamma^2)/2$;
\item
for $x\in (0,1/2)$, we compute $S(x)$ as in \eqref{S-periodic};
\item
for $x\in (1/2,1)$, we compute $S(x)$ as in \eqref{S-new-trunc-series}.
\end{enumerate}
\end{Corollary}
The proof of Corollary \ref{S-computation} follows just collecting the
information coming from Theorem \ref{main-thm-S}, Corollary \ref{corollary-S}, 
equations \eqref{S-value-12} and \eqref{difference-S}.

Even if in our application we will always work with $x\in(0,1)$, we recall that,
for $x$ large, it might be useful to implement the Stirling-like formula proved
in Theorem 2.11 of Deninger \cite{Deninger1984} which gives an asymptotic expression 
for $R(x)$ and, \emph{a fortiori}, for $S(x)$.

Our second theorem is about the function $T$ defined in \eqref{T-def}. 
As for $S(x)$, the starting point is the following 
\begin{Theorem}
\label{main-thm-T}
Let $x\in (0,2)$. Using the notation introduced in Theorem \ref{main-thm-S}, we
have
\begin{equation}
\label{T-new-series}
T(x) = 
\sum_{k=2}^{+\infty}\ \big( \zeta(k) H_{k-1} +\zeta^\prime(k)\bigl)(1-x)^{k-1}.
\end{equation}
Moreover, letting $x\in (0,1)\cup(1,2)$, $n\in \N$, $n\ge 1$ be fixed, and 
$r_T (x,n)\in \N$, 
\begin{equation}
\label{r-estim-T}
r_T (x,n):=
\Bigl\lceil
\min_{r}
\Bigl\{
r\ge 
1+ \frac{(n+2) \log 2 - \log \vert \log \vert 1-x \vert \vert+ 
\log \log r}{\vert \log\vert 1-x \vert \vert }  
\Bigr\}
\Bigr\rceil,
\end{equation}
we have that there exists
$\theta=\theta(x)\in (-1,1)$ such that
\begin{equation}
\label{T-new-trunc-series}
T(x) =
\sum_{k=2}^{r_T (x,n)} \big( \zeta(k) H_{k-1} +\zeta^\prime(k)\bigl)(1-x)^{k-1}
\ +\
\vert \theta\vert 2^{-n}.
\end{equation}
\end{Theorem}

We immediately remark that \eqref{T-new-series} is the Taylor series centred at $1$
of $T(x)$; in particular this implies  that
\[ 
T^{(k)} (1)=
 (-1)^k k! \bigl( \zeta(k+1) H_{k} + \zeta^\prime(k+1)\bigl)
\quad (k\in \N, k\ge 1).
\] 

Using Theorem \ref{main-thm-T} and \eqref{T-def}, for $x\in(0,2)$ we trivially have 
\[
\psi_1(x) =
-\gamma_1 
+
\sum_{k=2}^{+\infty}\ \big( \zeta(k) H_{k-1} +\zeta^\prime(k)\bigl)(1-x)^{k-1}
\]
and the corresponding truncated version.

Formula \eqref{T-new-series} is  essentially the one in Entry $21$(ii) on page $280$ of \cite{Berndt1985}
and it follows by differentiation
from \eqref{S-new-series} of Theorem \ref{main-thm-S}. The series in 
Theorem \ref{main-thm-T} clearly has a worst convergence speed than the one in Theorem \ref{main-thm-S}
and this justifies the 
different bound on $r_T (x,n)$ we have in \eqref{r-estim-T} comparing
with the one for $r_S (x,n)$ in \eqref{r-estim-S}. 
Recalling \eqref{TS-link-difference}, the fact that Theorem \ref{main-thm-T}
holds for every $x\in(0,2)$ means that every value of $T(x)$, $x\in(0,1)$, can be
computed in two different ways and that the shifting trick can be used in this case too.\footnote{We remark
that the size of the convergence interval can be doubled by isolating
the Taylor series at $1$ of $2(\log x)/x$ and using the estimates
on $\vert \zeta(n)-1 \vert$ of Lemma \ref{elementary-estim} below. We do not
insert such an idea here, since the computation of such an extra-factor $2(\log x)/x$ 
leads, in our practical application, to a longer total running time.} 
Hence,  if $x\in (1/2,1)$ we will directly compute $T(x)$ using Theorem
\ref{main-thm-T} while for $x\in (0,1/2)$ we will use
\eqref{TS-link-difference} and  Theorem \ref{main-thm-T} in $(1,3/2)$. 
This way we obtain the following two corollaries.
\begin{Corollary}
\label{corollary-T}
Let $x\in (0,1/2)$. We have that 
\begin{equation}
\label{T-periodic-series}
T(x) =
-
\frac{\log x}{x} 
+
\sum_{k=2}^{+\infty} (-1)^{k-1}\bigl(\zeta(k) H_{k-1} + \zeta^\prime(k)\bigl) x^{k-1}. 
\end{equation}
Letting further $n\in \N$, $n\ge 1$ be fixed
and  $r^\prime_T (x,n) \in \N$, 
\begin{equation*}
r^\prime_T (x,n)
:=
r_T (1+x,n) 
=
\Bigl\lceil
\min_{r}
\Bigl\{
r\ge 
1+ \frac{(n+2) \log 2 - \log \vert \log x \vert+ 
\log \log r}{\vert \log x \vert }  
\Bigr\}
\Bigr\rceil,
\end{equation*}
where $r_T(u,n)$ is defined in Theorem \ref{main-thm-T}, 
we have that there exists $\eta=\eta(x)\in (-1,1)$ such that
\begin{equation}
\label{T-periodic}
T(x) =
-
\frac{\log x}{x} 
+
\sum_{k=2}^{r^\prime_T (x,n)}  (-1)^{k-1}\bigl(\zeta(k) H_{k-1} + \zeta^\prime(k)\bigl) x^{k-1}\ 
+
\
\vert \eta\vert 2^{-n}.
\end{equation}
\end{Corollary}
Recalling equation (7.14) of Dilcher \cite{Dilcher1992} we also obtain
\begin{equation}
\label{T-value-12}
T\Bigl(\frac12\Bigr)
=
\gamma_1 + \psi_1\Bigl(\frac12\Bigr)  
=
 (\log 2)^2 + 2 \gamma \log 2.
\end{equation}
\begin{Remark}
In this case too we remark that 
evaluating twice $T(1/2)$ using Theorem \ref{main-thm-T}
(the first time directly and the second as $T(3/2)+2\log2$), we have, by summing
such formulae, that the summands having even indices vanish; thus we obtain
\[
\gamma 
= 
- \frac12 \log 2 + \frac12
+\frac{1}{2\log 2}
\sum_{\ell=1}^{n/2+4} \frac{ \zeta(2\ell+1) H_{2\ell} + \zeta^\prime(2\ell+1)
}{ 4^{\ell}} 
+ 
\tau 2^{-n+2},
\]
for some $\tau\in(-1/2,1/2)$.
Such a result
allows us to fast compute $\gamma$ with a precision of $n$ bits using about
$n/2$ steps. For example, using PARI/GP, we got $\gamma$ with a precision of
$1\,000$ decimal digits within 
$1$ minute and $10$ seconds 
of computation time on the Dell Optiplex machine previously mentioned.
In this case too there exist much faster algorithms to perform such a
computation, see again \cite{JohanssonB2019}.
\end{Remark}
Combining \eqref{T-new-trunc-series}, \eqref{T-periodic} and \eqref{T-value-12},
we obtain a very fast way of computing $T(x)$ for every 
$x>0$. We will see more about this in Section \ref{sect-implementation} but we 
also summarise the situation in the following
\begin{Corollary}
\label{T-computation}
We use the notation introduced in Theorem \ref{main-thm-T} and Corollaries
\ref{S-computation}-\ref{corollary-T}. 
We have:
\begin{enumerate}[i)]
\item
 $T(1)=T(2)=0$ and  $T(m) = \sum_{k=2}^{m-1} (\log k)/k$
for every $m\in \N$, $m\ge 3$;
\item
if $x>1$, $x\not \in \N$, we compute $T(x)$ as $T(x) = T(\{x\}) +
\sum_{k=0}^{\lfloor x \rfloor - 1} (\log (\{x\}+k)) / (\{x\}+k)$;
\item
$T(1/2) =  (\log 2)^2 + 2 \gamma \log 2$;
\item
if $x\in (0,1/2)$, we compute $T(x)$ as in \eqref{T-periodic};
\item
if $x\in (1/2,1)$,  we compute $T(x)$ as in \eqref{T-new-trunc-series}.
\end{enumerate}
\end{Corollary}
The proof of Corollary \ref{T-computation} follows just collecting the
information coming from Theorem \ref{main-thm-T}, Corollary \ref{corollary-T}, 
equations \eqref{T-value-12} and \eqref{TS-link-difference}. 

We finally remark that the shifting trick applies to
any function which can be defined as the solution of a difference equation, like
$S$ and $T$, and that it can be expressed via a power series whose convergence
interval is twice as large than the step of the difference equation. 
Another classical example of such a phenomenon\footnote{We used it in 
\cite{Languasco2020} to numerically study Littlewood's bounds on 
$\vert L(1,\chi)\vert$.} 
is the pair of functions given by $\log \Gamma$ and $\psi=\Gamma^\prime/\Gamma$, for which the analogues 
of the formulae \eqref{S-new-series} and \eqref{T-new-series} were first 
proved by Euler, see, \emph{e.g.}, Section 3 of the beautiful survey
of Lagarias \cite{Lagarias2013}.
But this also holds for further generalisations of Euler's Gamma
function like the ones studied by Dilcher in \cite{Dilcher1994}; in fact
$S$ and $T$ are the first and easier cases of such generalisations. 

Here we are mainly interested in $S$ and $T$ because of the number-theoretic
applications concerning the logarithmic derivative at $1$ of 
Dirichlet $L$-functions, see Section \ref{sect-applications}. 
There we will examine how to compute such quantities in a fast way and the 
number-theoretic consequences we can infer from such data.

The paper is organised as follows: in Sections \ref{ThS-proof}-\ref{ThT-proof}
we will respectively prove Theorems \ref{main-thm-S}-\ref{main-thm-T}. In Section
\ref{sect-applications} we will describe the problems in which the use of 
$S(x)$ and $T(x)$ is relevant.
In Section \ref{sect-implementation} we will discuss the computational costs
and some of the implementation
features of the formulae in Theorems \ref{main-thm-S}-\ref{main-thm-T} 
with respect to the applications too.
Finally, Section \ref{tables} is dedicated  
to show some figures about the applications 
described in Section  \ref{sect-applications}. 

 \medskip 
\textbf{Acknowledgements}. 
 The calculations here described in Section \ref{sect-implementation}
were performed using the University of Padova Strategic Research Infrastructure 
Grant 2017: ``CAPRI: Calcolo ad Alte Prestazioni per la Ricerca e l'Innovazione'',
\url{http://capri.dei.unipd.it}.
The authors would also like to thank the anonymous referees for their remarks
and suggestions.

\mbox{}\vskip-1.2truecm
\section{Proof of Theorem \ref{main-thm-S}}
\label{ThS-proof}
 
We start with the following lemmas that might have some independent interest too.
\begin{Lemma}
\label{T-gammas}
Let $x>0$ be fixed, $T(x)$ be defined as in \eqref{T-def}
and $\gamma_1$  as in \eqref{gamma1-zetasecond-def}.
Moreover let $\gamma$ be the Euler-Mascheroni constant.
Then we have
\begin{equation} 
T(x) 
\label{T-integral} 
=\gamma_1  - \int_{0}^{+\infty} 
\Bigl( e^{-t} - \frac{te^{-xt} }{1-e^{-t}} \Bigr)\frac{ \gamma+ \log t}{t} 
\ \dx t
\end{equation}
and
\begin{equation}
\label{int-gamma1}
\int_{0}^{+\infty} 
\Bigl(  e^{-t} - \frac{te^{-t} }{1-e^{-t}} \Bigr)\frac{\gamma+ \log t}{t} \ \dx
t
= \gamma_1 .
\end{equation}
\end{Lemma}
\begin{Proof} 
Using \eqref{S-alt-def} and \eqref{TS-link-difference}, a differentiation
immediately 
gives \eqref{T-integral}. The second part follows from the first 
using $T(1)=0$. 
\end{Proof}

\begin{Lemma}
\label{S-gammas}
Let $x>0$ be fixed, $S(x)$ be defined as in \eqref{S-def}, 
$T(x)$ be defined as in \eqref{T-def}
and $\gamma_1$ as in \eqref{gamma1-zetasecond-def}.
Let moreover $\gamma$ be the Euler-Mascheroni constant.
Then we have
\begin{equation}
\label{T-alt-def2}
T(x) 
= 
 \int_{0}^{+\infty} 
\bigl( e^{(1-x)t} - 1\bigr)\frac{\gamma+\log t}{e^t-1} \ \dx t
\end{equation}
and
\begin{equation}
\label{S-alt-def2}
S(x) 
= 
-2 (1-x) \gamma_1 
+
2 \int_{0}^{+\infty} 
\bigl( e^{(1-x)t} - 1-(1-x)t\bigr)\frac{\gamma+\log t}{t(e^t-1)} \ \dx t.
\end{equation} 
\end{Lemma}
\begin{Proof} 
Inserting \eqref{int-gamma1} into \eqref{T-integral} and performing a trivial 
computation on absolutely convergent integrals gives
\eqref{T-alt-def2}. Moreover,
an algebraic manipulation on equation \eqref{S-alt-def} 
immediately give 
\begin{align*}
S(x) &= 
2\int_{0}^{+\infty} 
(x-1) \Bigl(  e^{-t} - \frac{te^{-t} }{1-e^{-t}} \Bigr)\frac{\gamma+ \log t}{t}
\ \dx t
\\&
\hskip1cm
+
2\int_{0}^{+\infty} 
\Bigl(\frac{(x-1)te^{-t} +e^{-xt} - e^{-t}}{1-e^{-t}} \Bigr)\frac{\gamma+\log
t}{t} \ \dx t
\end{align*}
which is allowed since both integrals absolutely converge.
Recalling Lemma \ref{T-gammas} we see that the 
first integral is equal to $2(x-1)\gamma_1$.
Another algebraic manipulation on the second integral
proves \eqref{S-alt-def2}. We also remark that \eqref{T-alt-def2}
can also be obtained by differentiation from \eqref{S-alt-def2}
and \eqref{TS-link-difference}.
\end{Proof}

We will also need the following elementary estimates.
\begin{Lemma}
\label{elementary-estim}
Let $\gamma $ be the Euler-Mascheroni constant, $\psi(s)=
\Gamma^\prime/\Gamma(s)$ be the digamma function and
let $x>0$. Then 
\[
\log x - \frac{1}{x}< \psi(x) < \log x.
\]
Moreover, for every $k\in \N$, $k\ge 3$, we have
$\psi(k) + \gamma = H_{k-1}$,
\(
1+ 2^{-k}<\zeta(k)<1+2^{1-k}
\)
and, for $k\in \N$, $k\ge 4$, also that
\[
-\frac{\log2 +(2/3)\log 3}{2^k} < \zeta^\prime(k) <-\frac{\log 2}{2^k}.
\]
\end{Lemma}
\begin{Proof}
The first inequality follows from Theorem 5 of Gordon \cite{Gordon1994}.
The second part follows from \eqref{harmonic-def} and the fact that
$\psi(1)=-\gamma$ and $\psi(x+1) = \psi(x)+ 1/x$; hence
$\psi(k) + \gamma = \sum_{j=1}^{k-1} 1/j $
for every $k\in \N$, $k\ge 2$.
The estimate on $\zeta(k)$, $k\ge 3$, follows immediately from the definition
of the Riemann zeta-function in  $\Re(s) >1$, and the integral test.
Recalling that
$\zeta^\prime(s) = -\sum_{n=2}^{+\infty} (\log n)  n^{-s}$,  $\Re(s) >1$,
we have $-\zeta^\prime(k) > (\log 2) 2^{-k}$ for every $k\in \N $, $k\ge 2$. Moreover, using
 that $(\log x)/x$ is a decreasing sequence for every $x\ge e$,
the last part of the lemma  follows by  remarking
\[
-\zeta^\prime(k) =  \frac{\log 2}{2^k} + \sum_{n=3}^{+\infty} \frac{\log n}{ n^{k}}
<
 \frac{\log 2}{2^k} +  \frac{\log 3}{3} \sum_{n=3}^{+\infty} \frac{1}{ n^{k-1}}
 = 
  \frac{\log 2}{2^k} +  \frac{\log 3}{3} \Bigl(\zeta(k-1)-1- \frac{1}{2^{k-1}}\Bigr)
\]
and  using the inequality $\zeta(m) <1+2^{1-m}$, $m\in \N$, $m\ge 3$, previously proved.
\end{Proof}

The proof of Theorem \ref{main-thm-S} now starts from \eqref{S-alt-def2} of
Lemma \ref{S-gammas}.
Let $x\in (0,2)$ and, for every $t\in \R$, define $f(x, t) : = e^{(1-x)t} - 
1-(1-x)t$. Hence \eqref{S-alt-def2} becomes
\begin{equation}
\label{S-f-equiv}
S(x) 
=
-2 (1-x) \gamma_1 
+ 
2 \int_{0}^{+\infty} f(x,t)\frac{\gamma+\log t}{t(e^t-1)} \ \dx t.
\end{equation}
Writing the Taylor expansion at $0$ of $f(x, \cdot)$, we can easily get
that
\(
f(x,t) 
=
\sum_{k=2}^{+\infty} t^k (1-x)^{k}/k! 
\)
which holds for every $t\in \R$ and $x\in (0,2)$. 
Hence 
\begin{equation}
\label{int-series-exchange}
\int_{0}^{+\infty} f(x,t)\frac{\gamma+\log t}{t(e^t-1)} \ \dx t
= 
\sum_{k=2}^{+\infty} 
\frac{(1-x)^{k}}{k!}
\int_{0}^{+\infty} 
\frac{t^{k-1}(\gamma+ \log t)}{e^t-1} \ \dx t
\end{equation}
in which we exchanged the series and the integral signs by exploiting 
their absolute convergence.
Let now $s\in \C$, $\Re(s)>1$.
Recalling the classical formula 
\begin{equation}
\label{classical}
\int_0^{+\infty} \frac{t^{s-1}}{e^t-1} \ \dx t = \Gamma(s) \zeta(s),
\end{equation}
differentiating over $s$ we immediately get 
\begin{equation}
\label{int-form-zeta-diff}
\int_0^{+\infty} \frac{t^{s-1} \log t}{e^t-1} \ \dx t =
\Gamma^\prime(s) \zeta(s) + \Gamma(s) \zeta^\prime(s).
\end{equation}
Recalling that $\psi(s) =\Gamma^\prime/\Gamma(s)$, 
Lemma \ref{elementary-estim} and $\Gamma(k) =(k-1)!$,
by inserting \eqref{classical}-\eqref{int-form-zeta-diff} into
\eqref{int-series-exchange} we obtain
\begin{equation}
\label{zetas-psi}
\int_{0}^{+\infty} f(x,t)\frac{\gamma+ \log t}{t(e^t-1)} \ \dx t
=
\sum_{k=2}^{+\infty} 
\frac{\zeta(k) H_{k-1} + \zeta^\prime(k)}{k}(1-x)^{k} . 
\end{equation}
Hence \eqref{S-new-series} immediately follows by inserting \eqref{zetas-psi}
into \eqref{S-f-equiv}.
This proves the first part of Theorem \ref{main-thm-S}.
We now prove the second part of Theorem \ref{main-thm-S}.
From now on we denote
\begin{equation}
\label{L-def}
\LL(k):=\zeta(k)H_{k-1}+\zeta^\prime(k).
\end{equation}
Moreover letting $r\in \N$, $r\ge 3$, we define 
\[
\Sigma_S(r,x) := \sum_{k=2}^{r} 
\frac{\LL(k)}{k}(1-x)^{k} 
\]
and
\(
E_S(r,x) := 
\sum_{k=r+1}^{+\infty} 
(\LL(k)/k)(1-x)^{k}.
\)
Hence from \eqref{S-new-series} we get
\begin{equation}
\label{S-troncata}
S(x) =
-2 (1-x) \gamma_1 
+
2\Sigma_S(r,x) 
+
2E_S(r,x).
\end{equation}
Let now $n\ge 1$ be fixed. For every fixed $x\in (0,2)$, $x\ne 1$, we 
will find  $r=r_S (x,n)\in \N$ such that $\vert E_S(r,x)\vert < 0.25 \cdot
2^{-n}$.
Using Lemma \ref{elementary-estim} we obtain, for $k\ge r+1 \ge 4$, that
\[
\frac{\vert \LL(k)\vert}{k}
<
\frac{\zeta(k) H_{k-1}+ \vert \zeta^\prime(k) \vert}{k} 
< 
\frac{(1+2^{-r} )(\log (r+1)+\gamma)+0.72\cdot 2^{-r}}{r+1}
<
0.76
\]
and hence, using the well-known formula about the sum of a geometric
progression, we can write
\begin{equation}
\label{ES-estim}
\vert E_S(r,x)\vert
<
0.76
\sum_{k=r+1}^{+\infty} 
\vert 1-x \vert^{k}
=
0.76
\frac{\vert 1-x \vert^{r+1}}{1-\vert 1-x\vert}.
\end{equation}
We now look for $r =r_S (x,n) \in \N$, $r\ge 3$, such that
$\frac{\vert 1-x \vert^{r+1}}{1-\vert 1-x\vert} \le 2^{-n-2}$. An easy
computation reveals that
\begin{equation}
\label{r-condition-S}
r+1 \ge\frac{(n+2) \log 2 + \vert \log (1-\vert 1-x\vert)\vert}{ \vert \log 
\vert 1-x \vert \vert}
\end{equation}
suffices. 
The second part of Theorem \ref{main-thm-S} then follows
from \eqref{S-troncata}-\eqref{r-condition-S}.

\mbox{}\vskip-1.25truecm
\section{Proof of Theorem \ref{main-thm-T}}
\label{ThT-proof}
We already remarked that \eqref{T-new-series} follows via
\eqref{TS-link-difference}
from \eqref{S-new-series} but a direct proof can also be obtained starting from
\eqref{T-alt-def2}
and arguing as in the proof of Theorem \ref{main-thm-S}.
We now prove the second part of Theorem \ref{main-thm-T}.
Letting $r\in \N$, $r\ge 3$, and recalling \eqref{L-def}, we define 
\[
\Sigma_T(r,x) := \sum_{k=2}^{r} 
\LL(k)(1-x)^{k-1} 
\]
and
\(
E_T(r,x) := 
\sum_{k=r+1}^{+\infty} 
\LL(k)(1-x)^{k-1}.
\)
Hence from \eqref{T-new-series} we get
\begin{equation}
\label{T-troncata}
T(x)
= 
\Sigma_T(r,x) 
+
E_T(r,x).
\end{equation}
Let now $n\ge 1$ be fixed. For every fixed $x\in (0,2)$, $x\ne 1$, we 
will find  $r=r_T (x,n)\in \N$ such that $\vert E_T(r,x)\vert < 2^{-n}$.
Using Lemma \ref{elementary-estim} we have, for $k\ge r+1 \ge 4$, that 
\begin{equation}
\label{ET-estim1}
\vert E_T(r,x)\vert
<
\frac{2}{\vert 1-x \vert}
\sum_{k=r+1}^{+\infty} 
\vert 1-x \vert^{k} \log k. 
\end{equation}
Assuming that $r\ge 1/x$, we have that $\vert 1-x \vert^{k} \log k$ is a 
decreasing sequence for $k\ge r+1$;
hence a partial integration argument gives that
\begin{align}
\notag
\vert E_T(r,x)\vert
&
<
\frac{2}{\vert 1-x \vert}
\int_{r}^{+\infty} 
\vert 1-x \vert^{u} \log u\ \dx u
<
2 \frac{\vert 1-x \vert^{r-1}}{\vert \log \vert 1-x \vert\vert } 
\Bigl(\log r + \frac{1}{r\vert\log \vert 1-x \vert\vert} \Bigr)
\\&
\label{ET-estim2}
<
4 \frac{\vert 1-x \vert^{r-1}}{\vert\log \vert 1-x \vert\vert}\log r, 
\end{align}
in which we also assumed that $r\vert \log \vert 1-x \vert\vert\ge 1$.
We now look for $r =r_T (x,n) \in \N$ such that
$\frac{\vert 1-x \vert^{r-1}}{\vert\log \vert 1-x \vert\vert}\log r\le
2^{-n-2}$. 
An easy computation reveals that
\begin{equation}
\label{r-condition-T}
r-1 \ge\frac{(n+2) \log 2 - \log \vert \log \vert 1-x \vert \vert + \log \log
r}
{\vert \log\vert 1-x \vert \vert}
\end{equation}
suffices. The second part of Theorem \ref{main-thm-T} then follows
from \eqref{T-troncata}-\eqref{r-condition-T} . 

\mbox{}\vskip-1.25truecm
\section{Applications} 
\label{sect-applications} 
We briefly describe here some number-theoretic applications in which the 
use of $S$ and $T$ is relevant;
we will heavily refer to \cite{Languasco2021a} in which a more detailed 
presentation is given.

\subsection{Computation of $L^\prime/L(1,\chi)$}
\label{derlog-comput}
The main application in which is important to know the values of $S(a/q)$, where
$q$ is an odd prime and $a=1,\dotsc,q-1$,  is
to evaluate  the logarithmic derivative  at $1$ of the Dirichlet $L$-functions. 
Following the argument of Section 3 of \cite{Languasco2021a} we have,
for $\chi$ primitive and odd, that
\begin{align} 
\frac{L^\prime}{L} (1,\chi)
\label{odd-2}
&= 
\gamma 
+
\log ( 2 \pi ) 
+ 
\frac{1}{B_{1,\overline{\chi}} }
\sum_{a=1}^{q-1} 
\overline{\chi}(a)\log\Bigl(\Gamma\bigl(\frac{a}{q}\bigr)\Bigr)
\end{align}
and, for $\chi$ primitive and even, $\chi \ne \chi_0$, that
\begin{align} 
\frac{L^\prime}{L} (1,\chi)
\label{even-simplified}
&= 
\gamma + \log(2\pi)
-\frac{1}{2} 
\frac{\sum_{a=1}^{q-1} \overline{\chi}(a) \ S(a/q)}
{\sum_{a=1}^{q-1} \overline{\chi}(a)\log\bigl(\Gamma(a/q)\bigr)},
\end{align}
where
\(
B_{1,\chi}
: = 
(\sum_{a=1}^{q-1} a \chi(a) )/q
\)
is the first $\chi$-Bernoulli number,
$L(s,\chi)$ denote the Dirichlet $L$-functions, $\chi$ run over the non-principal
Dirichlet 
characters mod $q$ 
and $\chi_0$ is the principal Dirichlet character mod $q$.
Using the values of $T(a/q)$ we can alternatively write,
for every  $\chi\ne \chi_0$, that
\begin{equation}
\label{general-simplified} 
\frac{L^\prime}{L} (1,\chi)
=
- \log q
- 
\frac{\sum_{a=1}^{q-1} \chi(a) \ T(a/q)}
{\sum_{a=1}^{q-1} \chi(a)\ 
\psi (a/q) 
}.
\end{equation}
Ford, Luca and Moree \cite{FordLM2014} were the  first to use  \eqref{general-simplified} 
and the Fast Fourier Transform (FFT) method  to compute $\G_q$, see \S\ref{EK-constants-defs} for its definition.  
To see how the  summations over $a$ can be efficiently performed using the FFT method, see, \emph{e.g.}, Section 4 of \cite{Languasco2021a}, 
and for this approach  the use  of \eqref{odd-2}-\eqref{even-simplified} 
leads to a faster algorithm than the  one which uses \eqref{general-simplified}
because in the former case a \emph{decimation in frequency} strategy can be
applied. 
This essentially means that just the $(q-1)/2$ values of $S(a/q)+S(1-a/q)$ are
needed to perform the summation over $a$ in \eqref{even-simplified}; combining 
this with the use of suitable reflection formulae for $S(x)$ lead to gain a factor $4$
in the computational cost of generating the $S$-values with respect to the cost 
of generating the $T$-values, see also subsection \ref{comput-gain-DIF}.

It is a well-known fact, see, \emph{e.g.}, Corollary 10.18 of Montgomery-Vaughan
\cite{MontgomeryV2007}, that  the 
logarithmic derivative at $1$ of the Dirichlet $L$-functions
is connected with the  distribution of non-trivial 
zeros of $L(s,\chi)$. 
In literature, see, \emph{e.g.} Ihara-Murty-Shimura 
\cite{IharaMS2009}, it is interesting
to study the  extremal values 
of the logarithmic derivative at $1$ of the Dirichlet $L$-functions
under the assumption
of the Generalised Riemann Hypothesis. 

\subsection{Extremal values of $L^\prime/L(1,\chi)$}
For every odd prime $q$ we define
\begin{equation*} 
M^{\textrm{odd}}_q
:=\max_{\chi\, \textrm{odd}}\
\Bigl\vert\frac{L^\prime}{L} (1,\chi)\Bigr\vert \ ,
\quad 
M^{\textrm{even}}_q
:=\max_{\substack{\chi \neq \chi_0\\ \chi\, \textrm{even}}} \
\Bigl\vert\frac{L^\prime}{L} (1,\chi)\Bigr\vert\ , 
\quad
M_q=\max_{\chi \neq \chi_0} \
\Bigl\vert \frac{L^\prime}{L} (1,\chi)\Bigr\vert\ .
\end{equation*}

Hence we can compute $M_q= \max(M^{\textrm{odd}}_q, M^{\textrm{even}}_q)$  
using  \eqref{odd-2}-\eqref{even-simplified}.
Numerical values for $M_q$ were obtained in \cite{Languasco2021a} for every odd
prime $q\le 10^6$. Such data are in agreement with the estimate proved by Ihara-Murty-Shimura 
\cite{IharaMS2009} (please remark that our $M_q$ is denoted as 
$Q_m$ there) since they proved that $M_q \leq (2+\odi{1}) \log \log q$ as 
$q$ tends to infinity, under the assumption of the Generalised Riemann Hypothesis. 
On the other hand, Lamzouri, in a personal communication with
the first author, remarked that, by adapting the techniques in his paper \cite{Lamzouri2015}, 
one can show that
\(M_q\geq (1+\odi{1})\log\log q\)
for every sufficiently large prime $q$.

We will extend here the study of $M_q$  
to the larger interval $q\le \bound$; 
we can do so because of the much faster algorithm to compute $S(a/q)$ presented here.
A  similar study for $m_q:= \min_{\chi \neq \chi_0} \
\vert L^\prime/L (1,\chi)\vert\ $, is performed, both theoretically and numerically,
in Lamzouri-Languasco \cite{LamzouriL2020}.

\subsection{The Euler-Kronecker constants for prime cyclotomic fields}
\label{EK-constants-defs}
Let $q$ be an odd prime, $\zeta_q$ be a primitive $q$-root of unity,
$\zeta_{\Q(\zeta_q)}(s)$ be the Dedekind zeta-function of $\Q(\zeta_q)$. 
It is a well known fact that $\zeta_{\Q(\zeta_q)}(s)$ has a simple pole at $s=1$; 
writing the expansion of $\zeta_{\Q(\zeta_q)}(s)$ near $s=1$ as
\[
\zeta_{\Q(\zeta_q)} (s) = \frac{c_{-1}}{s-1} + c_0 + \Odi{s-1},
\]
the \emph{Euler-Kronecker constant of} $\Q(\zeta_q)$ is defined as  
\[
\lim_{s\to 1} \Bigl(\frac{\zeta_{\Q(\zeta_q)} (s)}{c_{-1}} - \frac{1}{s-1}\Bigr)=
\frac{c_0}{c_{-1}}.
\]
In this cyclotomic case 
we have that the Dedekind zeta-function can be written as
$\zeta_{\Q(\zeta_q)} (s)= \zeta(s) \prod_{\chi \neq \chi_0} L(s,\chi)$,
where $\zeta(s)$ is the Riemann zeta-function.
By logarithmic differentiation, we immediately get that
the \emph{Euler-Kronecker constant
for the prime cyclotomic field $\Q(\zeta_q)$} is
\begin{equation*}
\G_{q}
: =
\gamma
+
\sum_{\chi \neq \chi_0} \frac{L^\prime}{L} (1,\chi).
\end{equation*} 
Sometimes the quantity $\G_q$ is denoted as $\gamma_q$
but this conflicts with notations used in literature.
Another interesting quantity related to $\G_q$ is
the Euler-Kronecker constant $\G_q^+$ for $\Q(\zeta_q+\zeta_q^{-1})$, the
maximal real subfield of $\Q(\zeta_q)$. 
According to eq.~(10) of Moree \cite{Moree2018} it is defined as
\begin{equation*}
\G_q^+ : = \gamma 
+
\sum_{\substack{\chi \neq \chi_0\\ \chi\, \textrm{even}}} 
\frac{L^\prime}{L} (1,\chi).
\end{equation*}
An extensive study about the properties of $\G_{q}$ and $\G^+_{q}$ was 
started by Ihara \cite{Ihara2006,Ihara2008} and carried over by many others; 
we just recall here the papers by Ford-Luca-Moree \cite{FordLM2014} and Languasco
\cite{Languasco2021a}  because they both have some computational results on 
$\G_{q}$ and $\G^+_{q}$.
 
For both $\G_{q}$ and $\G^+_{q}$ it is interesting to find their negative
values since Ihara conjectured that both these quantities should be positive. 
Such a conjecture for $\G_{q}$ was disproved by Ford-Luca-Moree \cite{FordLM2014} 
(other two occurrences of $\G_{q}<0$ were detected in \cite{Languasco2021a}).
No negative values of $\G^+_{q}$ are known so far.
We will extend here the search for negative values of $\G_{q}$ and $\G^+_{q}$
to the large bound $q\le \bound$; 
in this way we also prove that there are no negative values for both 
$\G_q$ and $\G_q^+$ for every odd prime $q$ up to $\bound$. 
We also evaluate such quantities for some very large $q$.
We can do so because of the much faster algorithm to compute $S(a/q)$ 
presented here. 

\mbox{}\vskip-1.2truecm
\section{Implementation}
\label{sect-implementation}

We discuss here some implementation features of the formulae in
Corollaries \ref{S-computation} and \ref{T-computation}.
Since for $x>1$ we can reduce the problem of evaluating
$S(x)$, or $T(x)$, to a sum of a finite number of $\log$-values
plus $S(\{x\})$, or $T(\{x\})$, in this Section we assume that $x\in (0,1)$.

\subsection{Number of summands}
\label{number-summands}
We already remarked in the Introduction that,
from the estimates on $r_S (x,n)$ and $r_T (x,n)$ in Theorems
\ref{main-thm-S}-\ref{main-thm-T},
the number of needed terms we have to consider to have a $n$-bits digit
precision result becomes arbitrarily large as $x\to0^+$. To avoid this problem 
we can in practice use the  formulae in Corollaries \ref{S-computation} 
and \ref{T-computation}. 
In both corollaries it is clear that the worst cases for $r_S (x,n)$ and
$r^\prime_S (x,n) $ (and, respectively, for $r_T (x,n)$ and $r^\prime_T (x,n)$) 
are obtained when $x$ approaches $1/2$. 
Hence we can get any value of $S(x)$, $x\in (0,1)$,
with a precision of $n$ binary digits, with at most $n+2$ summands (assuming
that the needed $\log$ and $\LL(k)$ values can be obtained with the
same precision).
Analogously we can get any value of $T(x)$, $x\in (0,1)$,
with a precision of $n$ binary digits, with at most $n+4+\log\log(n+4)$ summands
(assuming that the needed $\log$ and $\LL(k)$ values can be obtained with the
same precision). 

\subsection{Precomputed coefficients}
In \eqref{S-new-trunc-series} and \eqref{T-new-trunc-series} we have a power
series whose coefficients involve  the values $\LL(k)$, 
$k\in \N$, $k\ge 2$ (see \eqref{L-def} for the definition of $\LL(k)$).
Hence in both cases such values can be precomputed, stored and reused for any
$x\in (0,1)$.
Moreover, the estimates in Lemma \ref{elementary-estim} imply
that $\vert \zeta(k) - 1 \vert<10^{-200}$ for $k\ge 160$
and $\vert \zeta^\prime(k) \vert <10^{-200}$ for $k\ge 420$. Hence, after about 
$420$ terms just the contribution of $H_{k-1}$ matters in 
\eqref{S-new-trunc-series} and \eqref{T-new-trunc-series}.
So, after few hundreds terms, the problem of obtaining $\LL(k)$ reduces
to being able to evaluate $H_{k-1}$.
We also remark that
 the computation of the needed first hundreds values of $\zeta(k)$ and 
$\zeta^\prime(k)$ can be performed, for instance, using PARI/GP.

Another nice aspect we have in \eqref{S-new-trunc-series} and
\eqref{T-new-trunc-series} is that the powers $(1-x)^k$ can be computed by 
recurrence, starting from $(1-x)^2$ and $1-x$, respectively.
The same clearly holds for equations \eqref{S-periodic} and \eqref{T-periodic}
too.

All these remarks also reveal  that the tasks of evaluating $S(x)$ and $T(x)$ are 
essentially as difficult as evaluating  $(\log x)^2$ and $(\log x)/x$, when $x$ is close to $0$.

\subsection{Reflection formulae for $S(x)$}
\label{comput-gain-DIF}
As mentioned in Section \ref{derlog-comput} and extensively explained in
Section 4 of \cite{Languasco2021a}, the use of the FFT algorithm is important 
to efficiently compute $\G_{q}$, $\G^+_{q}$ and $M_q$.
In particular, using $S(x)$, a decimation in frequency strategy can be
implemented and hence it is important to have the following 
\emph{reflection formulae} for $S$.

We directly express such formulae using Theorem \ref{main-thm-S}
and Corollary \ref{corollary-S}, or Corollary \ref{S-computation}, even if 
similar ones which use \eqref{S-def} and \eqref{S-alt-def} are also available
(such formulae were in fact used in \cite{Languasco2021a}, see Section 4.2
there).
\begin{Proposition}
\label{DIF-formulae-S}
Let $x\in (0,1)$, $x\ne 1/2$, $n\in \N$, $n\ge 2$, $r_1(x,n) = \lceil\frac{(n+2)
\log 2 + \vert \log (1-x) \vert}{ \vert \log x \vert } -1\rceil/2$
and $r_2(x,n) = \lceil\frac{(n+2) \log 2 +\vert \log x \vert}{\vert \log (1-x)
\vert}-1\rceil/2$.
Using \eqref{L-def} and the notations of Theorem \ref{main-thm-S} and Corollary
\ref{S-computation},
we have that there exists $\theta=\theta(x)\in (-1/2,1/2)$ such
that
\begin{align}
S(x)+S(1-x)
&= 
(\log x)^2
+
2 \sum_{\ell=1}^{r_1} \frac{\LL(2\ell)}{\ell}x^{2\ell} 
\label{S-DIF-x<1/2}
+
  \vert \theta\vert 2^{-n},
\quad \bigl(0<x < \frac12 \bigr),
\end{align}
and
\begin{align}
S(x)+S(1-x)
&= 
(\log(1-x))^2
+
2 \sum_{\ell=1}^{r_2} \frac{\LL(2\ell)}{\ell} (1-x)^{2\ell} 
\label{S-DIF-x>1/2}
+
 \vert \theta\vert   2^{-n},
\quad \bigl(\frac12<x<1\bigr).
\end{align}
\end{Proposition}
\begin{Proof}
Assume that $0<x < 1/2$; in this case we compute $S(x)$ with \eqref{S-periodic-series}
and $S(1-x)$ with \eqref{S-new-series}. 
Since the series  absolutely converge, their sum is   the series
having as summands the sum of their coefficients. Arguing as in 
\eqref{ES-estim}, remarking that  
$r_1(x,n)= r_S(1-x,n) /2=  r^\prime_S(x,n)/2$ and recalling \eqref{L-def},
we immediately have that \eqref{S-DIF-x<1/2} holds since
the odd summands vanish.
Assume that $1/2<x<1$; in this case we compute $S(x)$ with
\eqref{S-new-series} and $S(1-x)$ with \eqref{S-periodic-series}. 
 Arguing as for $x\in (0,1/2)$,   remarking that  
$r_2(x,n)= r_S(x,n)/2 = r^\prime_S(1-x,n)/2$ and
recalling \eqref{L-def}, we immediately have that \eqref{S-DIF-x>1/2} holds since
the odd summands vanish. This completes the proof.
\end{Proof}

The corresponding  series for \eqref{S-DIF-x<1/2}-\eqref{S-DIF-x>1/2} 
are
\begin{align}
S(x)+S(1-x)
&= 
(\log x)^2
+
2 \sum_{\ell=1}^{+\infty} \frac{\LL(2\ell)}{\ell}x^{2\ell} 
\label{S-DIF-x<1/2-series} 
\quad \bigl(0<x < \frac12 \bigr),
 \\
S(x)+S(1-x)
&= 
(\log(1-x))^2
+
2 \sum_{\ell=1}^{+\infty} \frac{\LL(2\ell)}{\ell} (1-x)^{2\ell} 
\label{S-DIF-x>1/2-series} 
\quad \bigl(\frac12<x<1\bigr).
\end{align}
Using Lemma \ref{elementary-estim}, it is easy to prove that $\LL(k)>0$ for every $k\in \N$, $k\ge 2$,
and hence \eqref{S-DIF-x<1/2-series}-\eqref{S-DIF-x>1/2-series} and 
\eqref{S-value-12} prove that $S(x)+S(1-x) > 0$ for every $x\in(0,1)$.

We remark that in Proposition \ref{DIF-formulae-S} we have $r_2(x,n)
= r_1(1-x,n)$ for $x\in (0,1)$ and hence the right hand side of \eqref{S-DIF-x>1/2} 
can be obtained from the right hand side of \eqref{S-DIF-x<1/2} just replacing
any occurrence of $x$ with $1-x$ and vice versa.
Analogous formulae, involving just the odd summands, can also be obtained for
$S(x) -S(1-x)$ but we omit them since they have no use in the applications 
here considered.

The use of Proposition \ref{DIF-formulae-S} in our application is four times
faster than using \eqref{T-new-trunc-series} and \eqref{T-periodic-series} for the 
following reasons:
\begin{enumerate}[i)]
\item exploiting the decimation in frequency strategy we just need to evaluate 
\eqref{S-DIF-x<1/2}-\eqref{S-DIF-x>1/2} at $x=a/q$, for every $a = 1,\dotsc,
(q-1)/2$, while \eqref{T-new-trunc-series} and \eqref{T-periodic-series} need to be 
evaluated for every $a = 1,\dotsc, q-1$.
This improves the computational cost by a factor $2$;
\item the cancellation of the odd terms we have in
\eqref{S-DIF-x<1/2}-\eqref{S-DIF-x>1/2}
leads to another gain of a factor $2$ in the computational cost with respect
to \eqref{T-new-trunc-series} and \eqref{T-periodic-series}
since we just need to use half of the summands (the ones with even indices);
\item in \eqref{S-DIF-x<1/2}-\eqref{S-DIF-x>1/2} the values of
the Riemann $\zeta$-function at even integers are required and for them
we can use the well-known exact formulae involving the Bernoulli numbers $B_k$:
\(
\zeta(2\ell) = (-1)^{\ell+1} \frac{B_{2\ell}(2\pi)^{2\ell}} {2(2\ell)!},
\)
for every $\ell\in \N$, $\ell\ge 1$, where the Bernoulli numbers $B_k$ are
defined using the following series expansion:
\(
\frac{t}{e^t-1} =\sum_{k=0}^{+\infty} B_k\frac{t^k}{k!}
\), $ \vert t \vert <2\pi$, 
see, \emph{e.g.}, Cohen's book \cite[chapter 9]{Cohen2007}.
\end{enumerate}
As we said before, the use of Proposition \ref{DIF-formulae-S}, if possible, 
is particularly efficient.
To compare the practical running times of using 
\eqref{S-DIF-x<1/2}-\eqref{S-DIF-x>1/2} with previous implementations,
which used the series/integral definitions of $S(x)$, see \eqref{S-def} and
\eqref{S-alt-def}, we compared the two PARI/GP scripts used to obtain 
$S(a/q)+S(1-a/q)$ for every  $a=1,\dotsc,(q-1)/2$ when
$q=305741$, $6766811$, $212634221$.
The gain in speed is huge, and it seems to improve as $q$ becomes larger:
we observed that the use of \eqref{S-DIF-x<1/2}-\eqref{S-DIF-x>1/2} leads to a
computation time for $S(a/q)+S(1-a/q)$ (with a precision of 128 bits)
for every $a=1,\dotsc,(q-1)/2$ which is respectively about $405$, $829$, $1206$ times  
faster for the three primes mentioned before.

Further practical experiments confirmed such a computational time gain;
we will see more on this in the next subsections.

\subsection{Computational costs for the problems of Section
\ref{sect-applications}}

The applications described in Section \ref{sect-applications} require to
evaluate  \eqref{S-new-trunc-series} over $x=a/q$, $a=1,\dotsc,q-1$.
Using the estimates in subsection \ref{number-summands}  we have
$r_S(x,n),r^\prime_S(x,n) \le  n+2$ for
every $x\in (1/2,1)$ and, respectively, $x\in (0,1/2)$. Hence  the total cost 
of evaluating  \eqref{S-new-trunc-series} over $x=a/q$, for every $a=1,\dotsc,q-1$,
 is $\Odi{q n}$ floating point products, 
with a precision of $n$ binary digits, and $\Odi{q}$ evaluation
of the logarithm function at rational points less than $1/2$. 
Since the remaining part of the computations in our applications are three 
 Fourier Transforms of length $\le q-1$, see \cite[Table 1]{Languasco2021a},
 having a cost of $\Odi{q\log q}$ floating point products each, this
proves that the total computational cost of our applications 
is $\Odi{q (n +\log q)}$ floating point products, with a precision of $n$ binary 
digits.
In practice, since in such applications we can use a decimation in frequency
strategy, Proposition \ref{DIF-formulae-S} let us
directly evaluate $S(a/q)+S(1-a/q)$ for every  $a=1,\dotsc,(q-1)/2$ 
thus reducing of the cost of such a step by a
factor of at least $4$.
A similar asymptotic estimate $\Odi{q (n +\log q)}$ holds also using 
\eqref{T-new-trunc-series} in the applications but in this 
case we cannot use the decimation in frequency strategy, see again
\cite{Languasco2021a}; hence in practice such an algorithm has a total cost which 
is about four times larger than the one which uses the $S$-function. 
Anyway, we will need such a $T$-function implementation for being able to double check 
the results.

\subsection{Actual implementation of the $S$ and $T$ formulae}
Using the C programming language, we implemented the formulae of Proposition
\ref{DIF-formulae-S} since they are the ones needed for the applications of
Section \ref{sect-applications}.  
The summation is performed combining 
the ``pairwise summation'' \cite{Higham1993}  algorithm
with Kahan's \cite{Kahan1965} method  (the minimal block
for the pairwise summation algorithm is summed using 
Kahan's method)
to have a good compromise between precision, computational
cost and execution speed.
To write here a practical computation time, we remark that for $q= 50\, 040\,
955\,631$ such an implementation computed $S(a/q)+S(1-a/q)$ for every 
$a=1,\dotsc,(q-1)/2$ with a precision of 128 bits in about nine hours using a
single computing core of an HP 
machine equipped with  4 x Eight-Core Intel(R) Xeon(R) CPU E5-4640 0 @ 2.40GHz, 
and 256GB of RAM.
For comparison, in this case the expected running time of the implementation
used in \cite{Languasco2021a} would be about 3475 days on the same machine 
mentioned before (about 9350 times slower).

Clearly this huge improvement let us evaluate the quantities described in
Section \ref{sect-applications} for some really large prime numbers and also to 
extend their knowledge  for every odd prime up to $\bound$. 
Moreover, to be able to double check the results obtained with the
$S$-function, we analogously implemented the
formulae of Corollary \ref{T-computation}.

\subsection{FFT implementation and computational results}
\label{FFT-implement}
To implement the FFT method we used the FFTW \cite{FFTW} package
which is also able to handle very large cases via its \texttt{guru64} interface. 
Moreover, to be able to store the large arrays of data we produce to
initialise the input sequences  involved in the FFT and their  outputs,
we used the \texttt{mmap}
UNIX system call to map such arrays on the hard disk instead of storing
them on the RAM during the execution of the C-programs. Thus we were able to enlarge the range
of possible computations we can perform far beyond the size of the
available RAM memory.
But that was not enough to handle the large case we would have liked to evaluate:
$q= 50\,040\,955\,631= 2 \cdot 5 \cdot  5\,004\,095\,563 +1$. 
We have chosen this prime number because its evaluation using the function $v(q)$,
defined in the next paragraph, see \eqref{vq-def},  is ``large'' enough 
to let us think that it might be a good candidate to have $\G_{q}<0$
($v(50040955631)=1.2194\dotsc$); please see Section \ref{tables} for more about
$v(q)$ and its link with the negativity of $\G_{q}$. 
Moreover, by analysing the prime  factor structure
of the known examples for which  $\G_{q}$ is negative, namely
$q=964477901, 9109334831, 9854964401$, see \cite{FordLM2014}
and \cite{Languasco2021a}, we see that such primes $q$ have 
all a ``large'' prime in the factorisation of $q-1$:
$964\,477\,901 = 2\cdot 5 \cdot  9\,644\,779 +1$,
$9\,109\,334\,831 =2\cdot 5 \cdot  910\,933\,483 +1$,
$9\,854\,964\,401 =2^4\cdot 5^2 \cdot  197 \cdot 125\,063 +1$.
These two motivations are hence a strong suggestion about the
negativity of $\G_{50040955631}$, even if they are not sufficient to be certain of this.

The presence of a ``large'' prime factor in the   factorisation of $q-1$
leads in fact to another problem in using the so-called \texttt{plan-generation}  step of 
the FFTW package. The \texttt{plan-generation}  step of FFTW is a procedure
in which FFTW  self-decides how to combine  several FFT algorithms
to obtain their best combination to solve the particular instance
of the problem the user is interested in. This procedure also depends on  the 
prime factorisation of the length of the transform $N$: in our case $N=q-1$, or $N=(q-1)/2$. When  $N$ has at least
one ``large'' prime factor,  as in our case, the plan-generation step
might be very demanding in term of memory usage
(RAM). To overcome this, we have then to insert the use of \texttt{mmap} in the body
of the FFTW code to be able to divert the memory usage of the \texttt{plan-generation}  step
from the RAM to the hard disk. Clearly this increases the actual computation time but, at the same time, 
let us handle much larger cases, since, essentially, it is much easier, and cheaper, 
to retrieve large hard disks than a large quantity of RAM.

In this way we obtained a program that needed at most 128GB of RAM at runtime
and we used it to perform the computation for the case of $q= 50\, 040\, 955\,631$
on the University of Padova Strategic Research Infrastructure  
``CAPRI'' (Intel(R) Xeon(R) Gold 6130 CPU @ 2.10GHz, with 256 cores and equipped
with 6TB of RAM).
The total hard disk usage was about 8 TB, 
the time needed  for one computing core to generate the $S$-values with a precision of 128 bits was about six hours
and 6 minutes, the \texttt{plan-generation} step required about four hours and 5 minutes
and the actual FFT transforms  about 2 days and half (for S).
The total computation time was about two weeks; we recall that
such computation times are affected, as above remarked, from the fact that
we were using a slower memory device (the hard disk is used instead of RAM).
We got that 
$\G_{50 040  955 631} =-0.16595399\dotsc$ and
$\G_{50 040  955 631}^+ =13.89764738\dotsc$
thus getting another occurrence of a negative Euler-Kronecker constant.

The computation of $\G_{q}, \G_{q}^+$ and $M_q$ for every odd prime $q$ up to $\bound$ 
was performed on CAPRI using at most  60 computing nodes and 
it required about 48 hours of  time 
(the global execution time, obtained by summing the declared computing time on 
each node, was of 101 days and 6 hours).

In this range we obtained that there are no negative values for both $\G_{q}$ and $\G_{q}^+$
and that
\begin{equation}
\label{Mq-comput-bounds}
\frac{17}{20} \log \log q< M_q < \frac{5}{4} \log \log q
\end{equation}
for every  prime $1531 < q\le \bound$;
$M_{1531}=2.5048094\dots$, $M_{1531}/ \log \log(1531)=1.257133\dots$ 
Moreover the lower bound  in \eqref{Mq-comput-bounds}
holds true for  $q>13$. 
 
The programs used and the results here described are collected 
at the following address
\url{http://www.math.unipd.it/~languasc/Scomp-appl.html}.

\subsection{FFT accuracy estimate}
According to Schatzman \cite[\S~3.4, p.~1159-1160]{Schatzman1996},  
the root mean square relative error in the FFT is bounded by  
\begin{equation}
\label{Delta-FFT}
\Delta = \Delta(N, \eps) :=0.6 \eps (\log_2 N)^{1/2} ,
\end{equation} 
where $\eps$ is the machine epsilon
and $N$ is the length of the transform. 
According to the IEEE 754-2008 specification, we can 
set $\eps=2^{-64}$ for the \emph{long double precision} 
of the C programming language. 
So for the largest case we are considering, $q=50 040 955 631$, $N=(q-1)/2$, we get that $\Delta<1.92 \cdot 10^{-19}$.
To evaluate the euclidean norm of the error we have then to multiply $\Delta$ and the
euclidean norms of the  sequences (listed in sections 4.2-4.3 of \cite{Languasco2021a}):
\begin{alignat*}{3}
x_k&:= 2\frac{a_k}{q}-1, \quad
&y_k&:= \log\Gamma\Bigl(\frac{a_k}{q}\Bigr) + \log\Gamma\Bigl(1-\frac{a_k}{q}\Bigr) - \log \pi,\\
z_k&:= \log\Gamma\Bigl(\frac{a_k}{q}\Bigr) - \log\Gamma\Bigl(1-\frac{a_k}{q}\Bigr),\quad
&w_k &:=S\Bigl(\frac{a_k}{q}\Bigr) -S\Bigl(1-\frac{a_k}{q}\Bigr),
\end{alignat*}
where  $a_k  = g^k \bmod q$, $\langle q \rangle =\Z_q^*$, $k=0,\dotsc, N-1$.
A straightforward computation gives 
\[
\Vert x_k \Vert_2 = \Bigl(   \frac{(q-1)(q-2)}{6q} \Bigr)^{1/2}
=91324.47246\dotsc
\]
Hence, recalling that $\Vert \cdot \Vert_{\infty} \le \Vert \cdot \Vert_{2}$,
 for this sequence we can estimate that the maximal error in its FFT-computation
is bounded by $ 1.75 \cdot 10^{-14}$ (long double precision case).
Unfortunately, no closed formulas for the euclidean norms of the other involved sequences are known 
but, using  $\Vert \cdot \Vert_{\infty} \le \Vert \cdot \Vert_{2} \le \sqrt{N} \Vert \cdot \Vert_{\infty}$
and the formulae
\begin{align*}
\Vert y_k \Vert_{\infty} &=   
- \log \sin(\pi/q) = 
23.49137\dotsc,\\
\Vert z_k \Vert_{\infty} &= 
2  \log\Gamma\Bigl( \frac{1}{q} \Bigr) 
-
\log\Bigl(\frac{\pi}{\sin(\pi/q)}\Bigr)=
24.63610\dotsc,\\
\Vert w_k \Vert_{\infty} &
= S\Bigl(\frac1q\Bigr) + S\Bigl(1-\frac1q\Bigr)
= 606.93779\dotsc,
\end{align*}
that can be obtained using straightforward computations,
we have that 
the errors in their  FFT-computations are all
$< 1.85 \cdot 10^{-11}$.

We also estimated \emph{in practice} the accuracy in 
the actual computations  using the 
FFTW software library by evaluating at   
run-time the quantity
\(
\mathcal{E}_{j}(w_k) : = 
\Vert \mathcal{F}^{-1}(\mathcal{F}(w_k)) - w_k \Vert_{j},
\)
$j\in\{2,\infty\}$, 
$\mathcal{F}(\cdot)$ is the Fast Fourier Transform 
and $\mathcal{F}^{-1}(\cdot)$ is its inverse transform.
We focused our attention on $w_k$  since,  between 
the sequences mentioned before, it has the largest norms
and hence the worst error estimates.
Theoretically we have that  $\mathcal{E}_{j}(w_k) =0$;
moreover, assuming that the root mean square relative error in the FFT is bounded by  
$\Delta>0$, it is easy to obtain
\begin{equation}
\label{back-forth-estim-general}
\mathcal{E}_{2}(w_k) < \Delta(2+\Delta) \Vert w_k \Vert_{2}
\quad
\textrm{and}
\quad
\mathcal{E}_{\infty}(w_k) < \Delta(2+\Delta)  \sqrt{N}  \Vert w_k \Vert_{\infty}.
\end{equation}
For  $q=50 040 955 631$, $N=(q-1)/2$ and $\eps = 2^{-64}$ in \eqref{Delta-FFT},
we get $\Delta(2+\Delta)< 3.83 \cdot 10^{-19}$ and 
from \eqref{back-forth-estim-general} we obtain
\begin{equation}
\label{back-forth-estim-50040955631}
\mathcal{E}_{2}(w_k) <  3.70 \cdot 10^{-11}
\quad
\textrm{and}
\quad
\mathcal{E}_{\infty}(w_k) < 3.70 \cdot 10^{-11},
\end{equation}
where the first estimate suffers from the lack of theoretical 
information about $\Vert w_k \Vert_{2}$.
Moreover, the actual computations using FFTW for this case gave that 
$\Vert w_k \Vert_2 =1099611.166707\dotsc$,  
\begin{equation} 
\label{E_wk_computed}
\frac{\mathcal{E}_2(w_k)}{\Vert w_k \Vert_2} < 6.01\cdot 10^{-19},
\quad
\mathcal{E}_2(w_k)< 4.21\cdot 10^{-13}
\quad
\textrm{and}
\quad
\mathcal{E}_\infty(w_k) < 2.23 \cdot 10^{-16}
\end{equation}
that are  in agreement with \eqref{back-forth-estim-50040955631}.
It is worth notice that the last two computed estimates in \eqref{E_wk_computed}
are much better than the corresponding theoretical ones in \eqref{back-forth-estim-50040955631}.
We finally remark that the computed estimates for the analogous quantities involving $x_k,y_k,z_k$ 
are smaller than the ones for $w_k$ described before.

Summarising, we can conclude that at least ten decimal 
digits of our final results are correct.
If necessary, more accurate results can be obtained using the 
\emph{quadruple precision} ($128$ bits),
which let us choose $\eps = 2^{-113}$ in \eqref{Delta-FFT},
at the cost of a much slower execution. 

We finally recall
that the well-known weakness of the FFT-algorithms is 
the memory occupation and, as another example 
of this, we remark that for $q=50 040 955 631$ about 12 TB of 
hard disk space are required (see also \S \ref{FFT-implement})
 to perform the needed FFTW-library computations
(long double precision case) used to obtain the data described in
this paragraph.

\mbox{}\vskip-1truecm
\section{Figures}
\label{tables} 
We give here some comments about Figures \ref{fig1}-\ref{fig4} and discuss the role of the $v(q)$-function.
 
Referring to Section 4.5 of \cite{Languasco2021a}, we recall the definition
of   $\B$, the ``greedy sequence of prime offsets'', 
 \url{http://oeis.org/A135311}.
 We define $\B$ using induction, by 
   $b(1)=0\in \B$ and $b(n)\in \B$ if it is  the smallest integer exceeding $b(n-1)$ such that 
 for every prime $r$ the set $\{b(i) \bmod r\colon 1\le i \le n\}$ has at most $r-1$ elements. 
Let now
\[m(\A):= \sum_{i=1}^s\frac{1}{a_i},\]
 where $\A$ is an admissible set, \emph{i.e.}, $\A=\{a_1,\dotsc, a_s \}$, $a_i\in \N$, $a_i\ge 1$, such that
 does not exist a prime $p$ such that $p\mid n\prod_{i=1}^s (a_in+1)$ for every $n\ge 1$.
Thanks to Theorem 2 of
Moree \cite{Moree2018}, if the prime $k$-tuples conjecture holds and if $\A$ is an admissible set, then
$\G_q < (2- m(\A)+ \odi{1})  \log q$ for $\gg x/(\log x)^{-\vert \A \vert -1}$ primes $q \le x$.
Moreover, by Theorem 6 of Moree \cite{Moree2018}, assuming both the Elliott-Halberstam 
and the prime $k$-tuples conjectures,  if $\A$ is an admissible set then
$\G_q = (1- m(\A)+ \odi{1})  \log q$ for $\gg x/(\log x)^{-\vert \A \vert -1}$ primes $q \le x$.
We recall that  the greedy sequence of prime offsets $\B$ has the property that any finite 
subsequence is an admissible set.
 With a PARI/GP script we   computed the first $2089$ 
elements of $\B$  since for $\CC := \{b(2), \dotsc, b(2089)\}$ we get  $m(\CC) >2$.
 
\begin{figure}[H]
 \begin{minipage}{0.48\textwidth}
 \includegraphics[scale=0.525,angle=0]{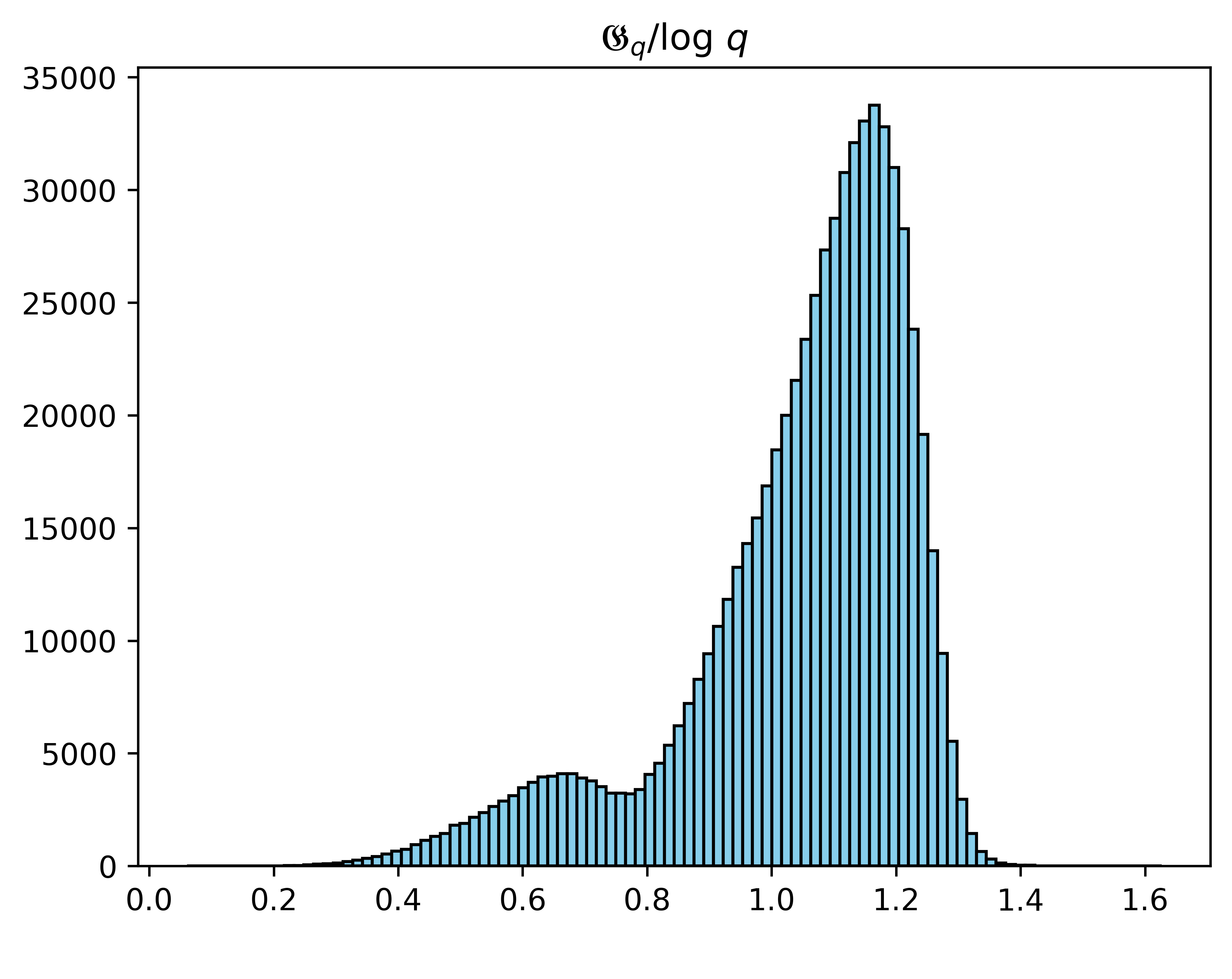}   
\end{minipage}
\begin{minipage}{0.48\textwidth} 
  \includegraphics[scale=0.525,angle=0]{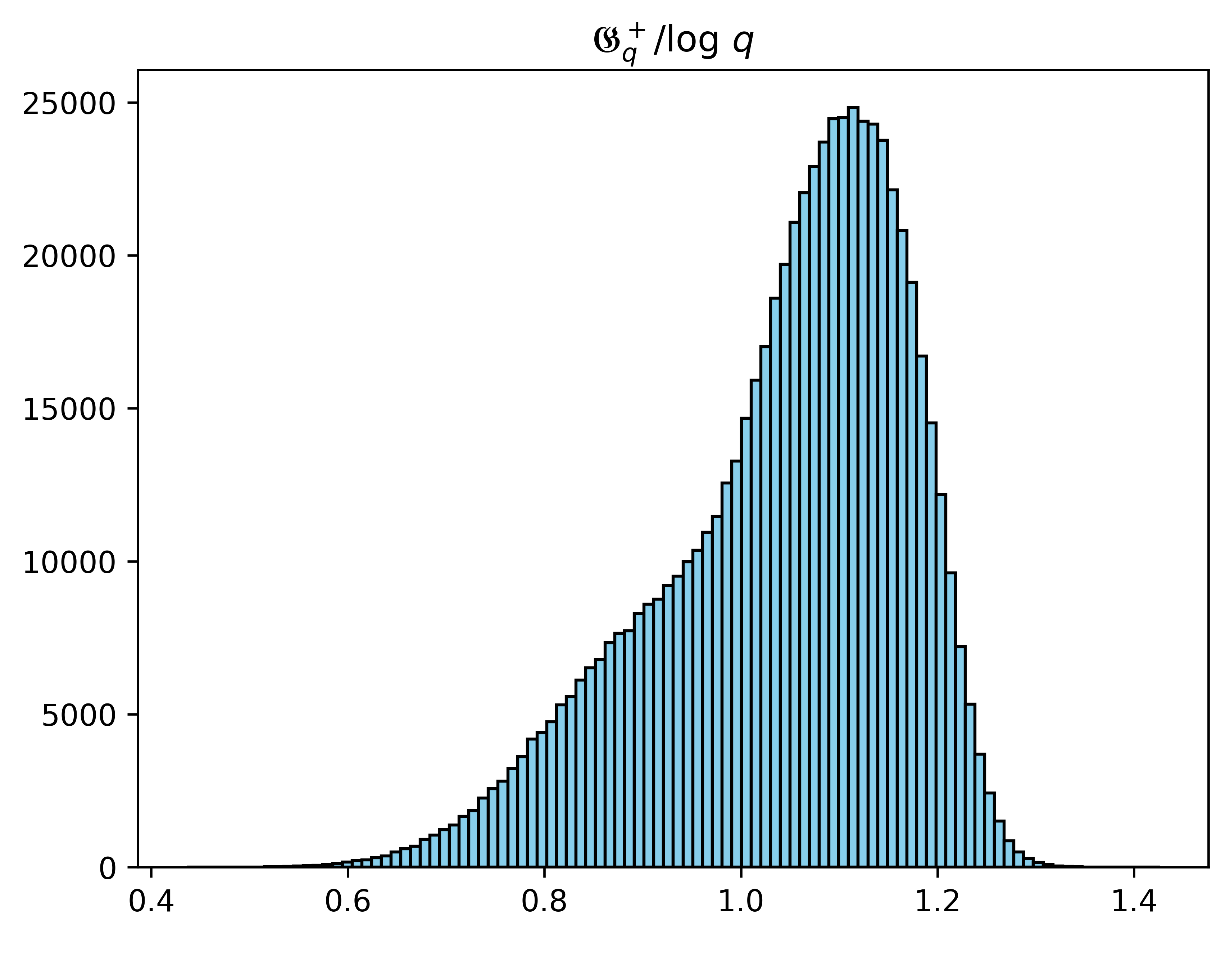}   
\end{minipage} 
\caption{\small{The histograms for $\G_{q}/\log q$ and $\G_{q}^+/\log q$,   $q$ prime, $3\le q\le  \bound$.
 In both the
number of intervals is $100$, the
number of odd primes up to $10^7$ is ${\mathcal P} = 664578$
and the mass is $\mathcal{M} =  I \cdot {\mathcal P}$,
where $I$ is the interval length.
For the first we have $I =  0.017230\dotsc$,
the mean is $\mu = 1.043452\dotsc$
and  the
standard deviation is $\sigma =0.185900\dotsc$
For the second we have $I =  0.010892\dotsc$,
the mean is $\mu = 1.043347\dotsc$
and  the
standard deviation is $\sigma =0.125566\dotsc$
}}
 \label{histoEK-EK+} 
 \end{figure}

\begin{figure} [H]  
\begin{minipage}{0.48\textwidth}
\centerline{{\tiny $\G_q/ \log q$ with $v(q)\le 0.25$}}
\includegraphics[scale=0.625,angle=0]{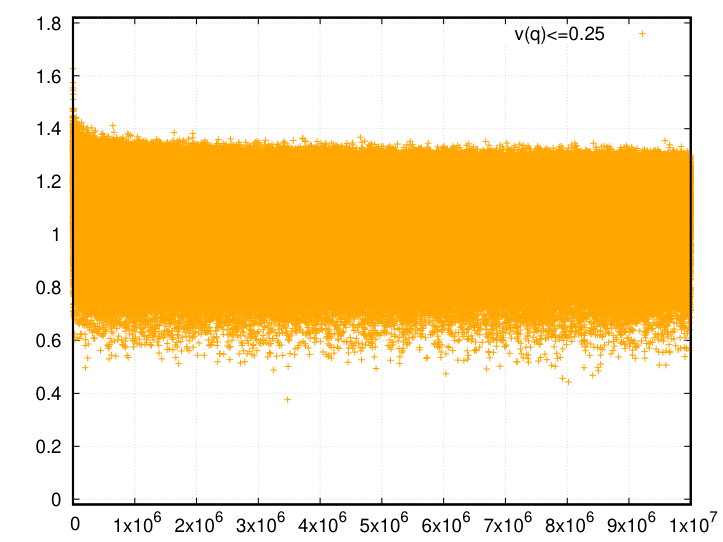}
\end{minipage}
\begin{minipage}{0.48\textwidth}
\centerline{{\tiny $\G_q/ \log q$ with $0.25 < v(q)\le 1/2$}}
\includegraphics[scale=0.625,angle=0]{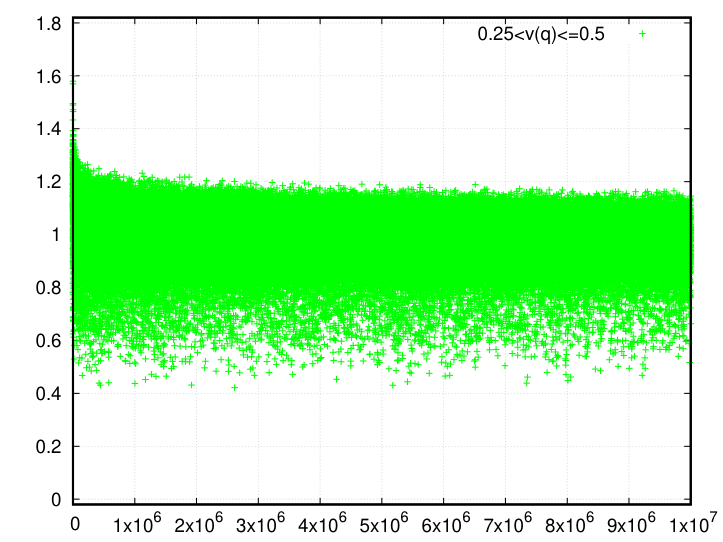}
\end{minipage} 

\vskip1.5cm
\begin{minipage}{0.48\textwidth}
\centerline{{\tiny $\G_q/ \log q$ with $1/2 < v(q)\le 0.75$}}
\includegraphics[scale=0.625,angle=0]{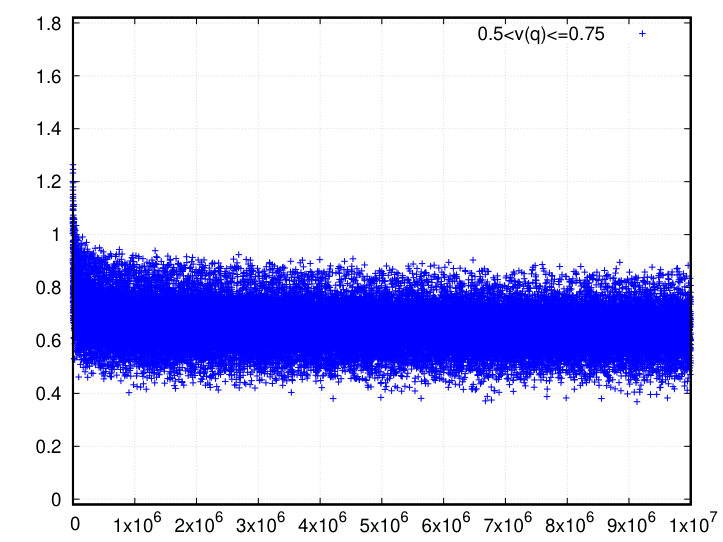}  
\end{minipage}
\begin{minipage}{0.48\textwidth}
\centerline{{\tiny $\G_q/ \log q$ with $0.75 < v(q)\le 1$}}
\includegraphics[scale=0.625,angle=0]{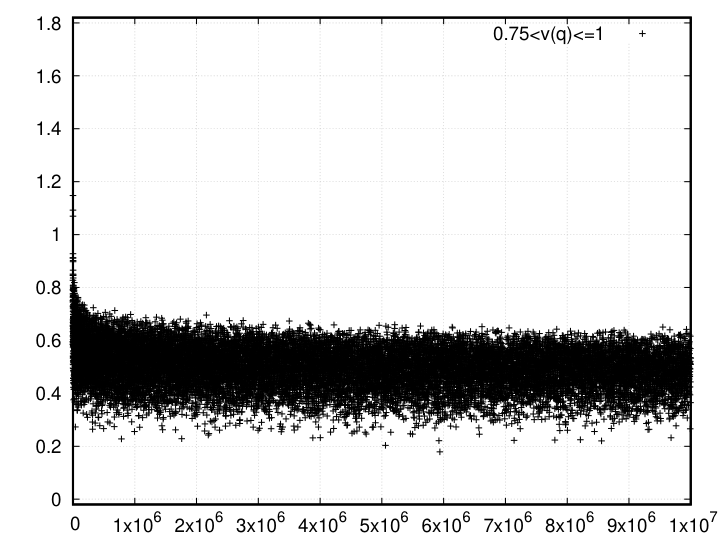}  
\end{minipage}

\vskip1.5cm
\begin{minipage}{0.48\textwidth}
\centerline{{\tiny $\G_q/ \log q$ with $v(q) > 1$}}
\includegraphics[scale=0.625,angle=0]{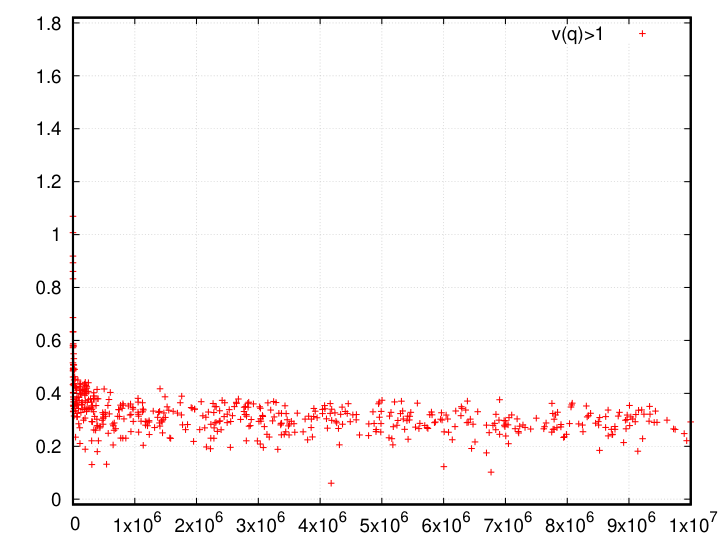}  
\end{minipage}
\begin{minipage}{0.48\textwidth}
\centerline{{\tiny $\G_q/ \log q$}}
\includegraphics[scale=0.625,angle=0]{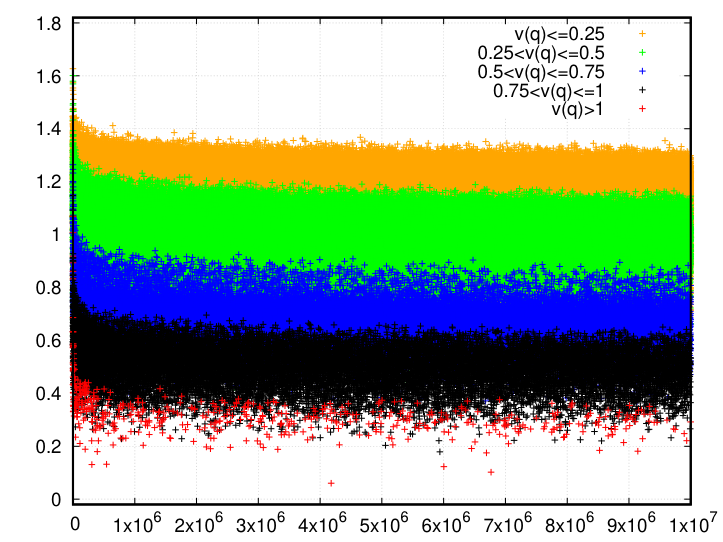}  
\end{minipage} 
\caption{{\small The values of $\G_q/\log q$, $q$ prime, $3\le q\le \bound$,
classified using $v(q)$.
The  minimal value is $0.060532\dotsc$ and it is attained at $q=4178771$; 
the  maximal value is $1.626934\dotsc$ and it is attained at $q=19$.}}
\label{fig1}
\end{figure}

\begin{figure} [H] 
\begin{minipage}{0.48\textwidth}
\centerline{{\tiny $\G^+_q/ \log q$ with $v(q)\le 0.25$}}
\includegraphics[scale=0.625,angle=0]{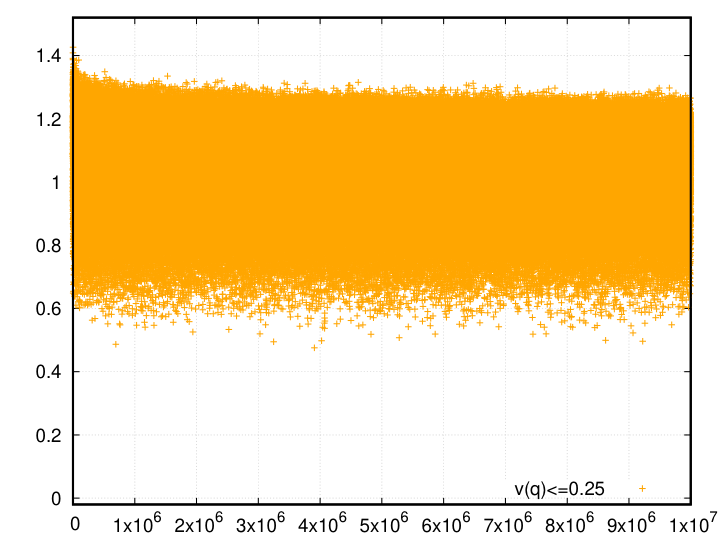} 
\end{minipage} 
\begin{minipage}{0.48\textwidth}
\centerline{{\tiny $\G^+_q/ \log q$ with $0.25 < v(q)\le 1/2$}}
\includegraphics[scale=0.625,angle=0]{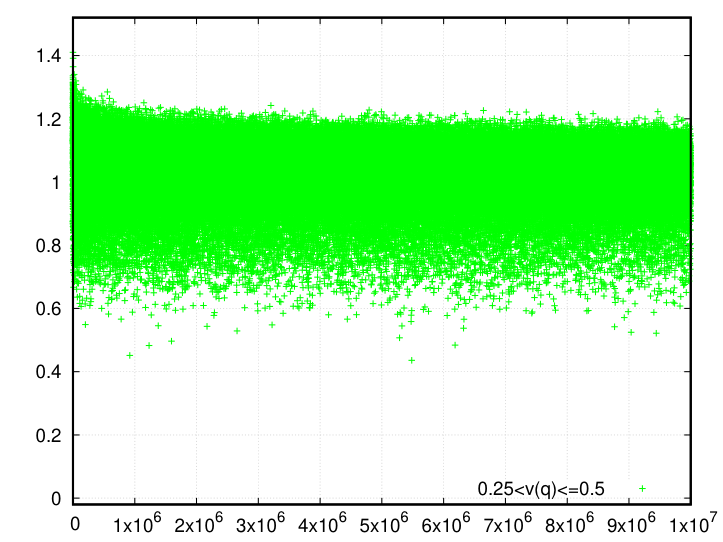}
\end{minipage} 

\vskip1.5cm
\begin{minipage}{0.48\textwidth} 
\centerline{{\tiny $\G^+_q/ \log q$ with $1/2 < v(q)\le 0.75$}}
\includegraphics[scale=0.625,angle=0]{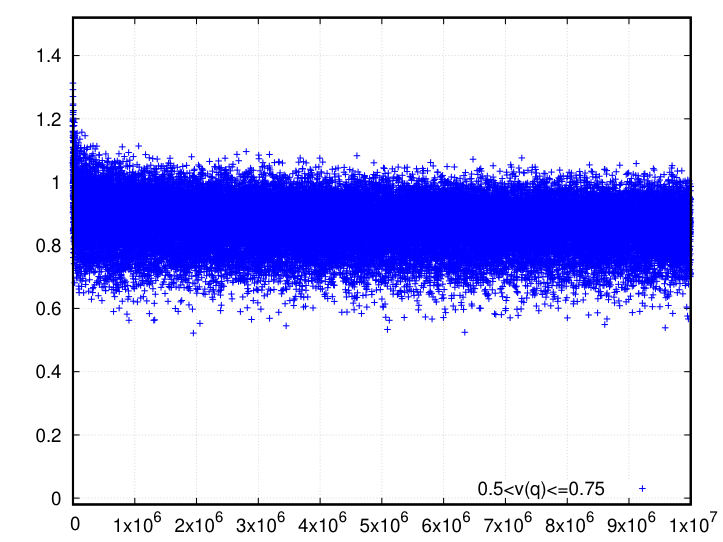}
\end{minipage} 
\begin{minipage}{0.48\textwidth} 
\centerline{{\tiny $\G^+_q/ \log q$ with $0.75 < v(q)\le 1$}}
\includegraphics[scale=0.625,angle=0]{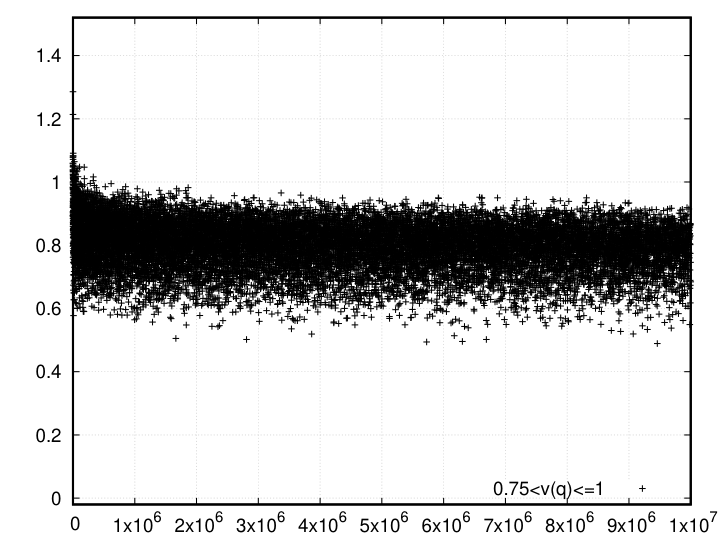} 
\end{minipage} 

\vskip1.5cm
\begin{minipage}{0.48\textwidth} 
\centerline{{\tiny $\G^+_q/ \log q$ with $v(q)>1$}}
\includegraphics[scale=0.625,angle=0]{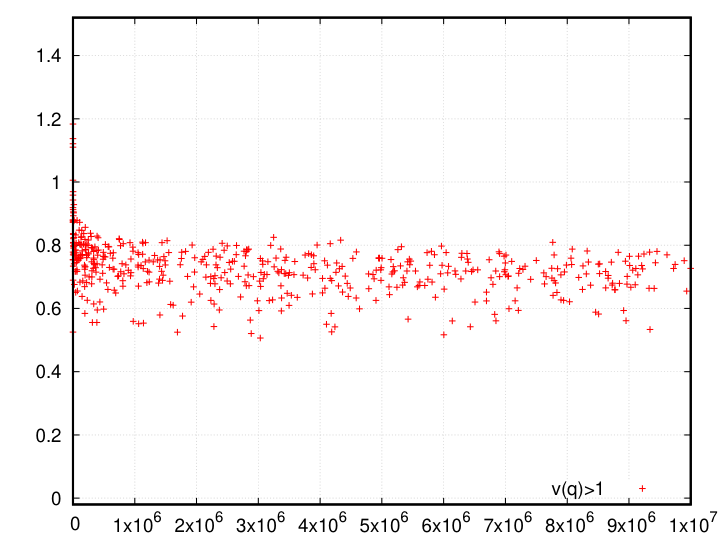} 
\end{minipage} 
\begin{minipage}{0.48\textwidth} 
\centerline{{\tiny $\G^+_q/ \log q$}}
\includegraphics[scale=0.625,angle=0]{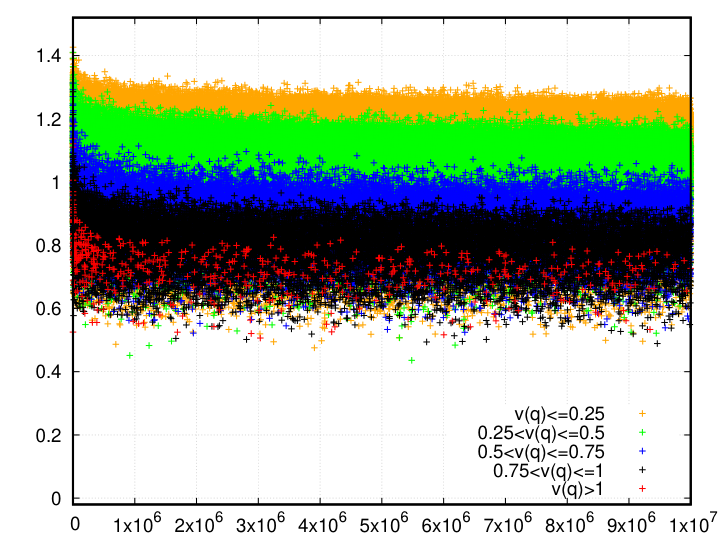}  
\end{minipage} 
 \caption{{\small The values of $\G_q^+/\log q$, $q$ prime, $3\le q\le \bound$,
 classified using $v(q)$.
 The  minimal value is $0.436031\dotsc$ and it is attained at $q= 5483977$; 
 the  maximal value is $1.426263\dotsc$ and it is attained at $q=2053$.}}  
 \label{fig2}
 \end{figure}

\begin{figure}[H] 

\begin{minipage}{0.48\textwidth} 
\centerline{{\tiny $M_q$, $3\le q\le 2000$}}
\includegraphics[scale=0.625,angle=0]{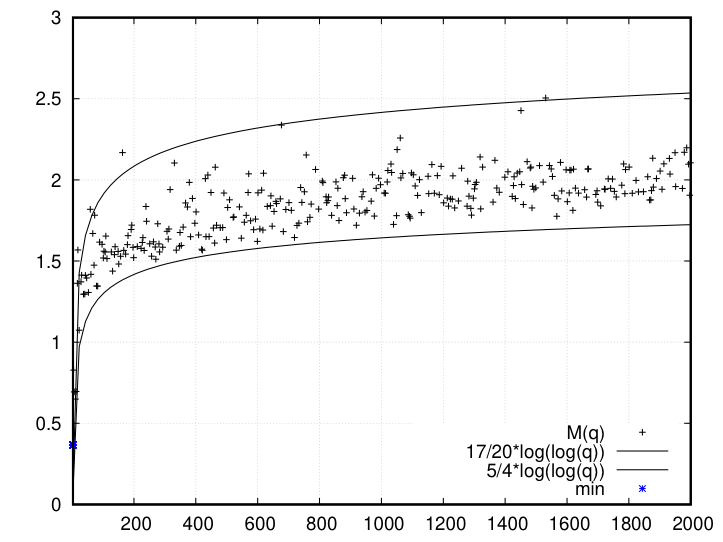}  
\caption{{\small The values of $M_q$, $q$ prime, $3\le q\le 2000$. 
 The  minimal value is $0.3682816 \dotsc$ and it is attained at $q=3$.
 The lines  represent  the functions $c\cdot\log \log q$, with $c=17/20$ and $c=5/4$.
 $M_{1531}=2.5048094\dots$
 }}
 \label{fig3a} 
\end{minipage} 
\begin{minipage}{0.48\textwidth} 
\centerline{{\tiny $M_q$, $2000< q\le 10^7$}}
\includegraphics[scale=0.625,angle=0]{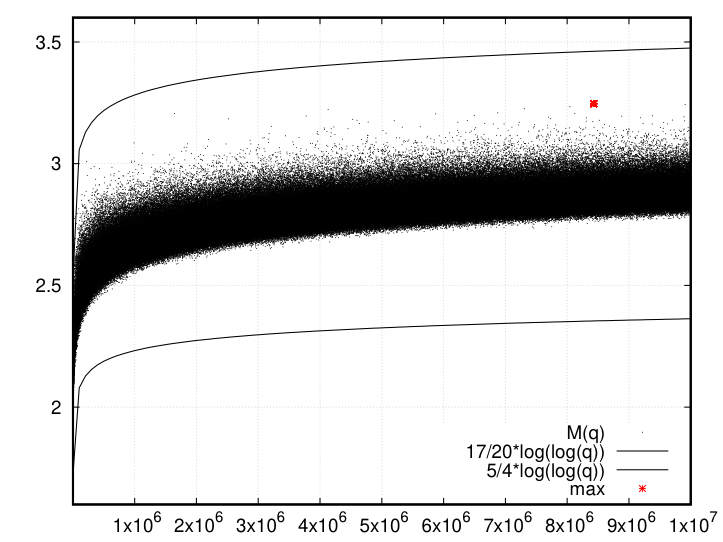}  
\caption{{\small The values of $M_q$, $q$ prime, $2000< q\le \bound$. 
The  maximal value is $3.2466918\dotsc$ and it is attained at $q= 8430391$. 
The lines  represent  the functions $c\cdot\log \log q$, with $c=17/20$ and $c=5/4$. 
}}
\label{fig3b} 
\end{minipage}  
\end{figure}

\begin{figure}[H] 
\centerline{{\tiny $M_q/\log \log q$}}
\includegraphics[scale=0.95,angle=0]{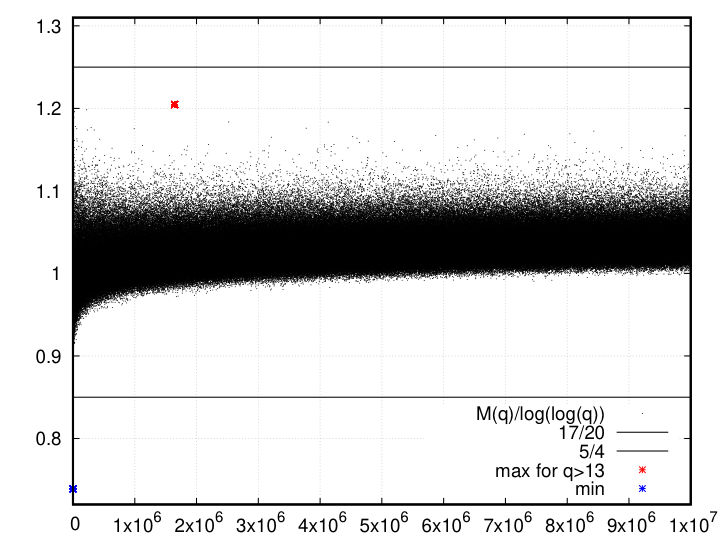}  
\caption{{\small The values of $M^\prime_q:=M_q/\log \log q$, $q$ prime, $3\le q\le \bound$. 
 The  minimal value is $0.7392305\dotsc$ and it is attained at $q=13$; 
 the  maximal value is $3.9158971\dotsc$ and it is attained at $q= 3$ (not represented
 in the plot); the maximal value for every $q>13$ is $1.204704\dotsc$
 and it is attained at $q= 1645093$.
 The  lines  represent  the constant functions  $c=17/20$ and $c=5/4$. 
 }}
\label{fig4} 
\end{figure}

So, if we are looking for negative values of $\G_q$, it seems to be a good criterion to evaluate $\G_q$ for 
a  prime number $q$  such that $bq+1$ is prime for many elements $b\in \CC$ 
(clearly it is better to start with the smaller available $b$'s).
 To be able to measure this  fact, we define
\begin{equation}
\label{vq-def}
v(q) :=  \sum_{\substack{2\le i\le 2089; \ b(i)\in \CC\\ b(i)q+1\ \text{is prime} }}\frac{1}{b(i)}
\end{equation}
and we use such a function to classify  $\G_q$ and $\G_q^+$ in the following way.
In the scatter plots of Figures \ref{fig1}-\ref{fig2} we classified the normalised values of $\G_q$ and $\G_q^+$
according to $v(q)$.   Orange points  are the more frequent ones ($72.88 \%$ of the total number) 
and satisfy $v(q)\le 0.25$; 
green points satisfy $0.25 <v(q)\le 1/2$ ($18.29\%$);    
blue points satisfy $1/2 <v(q)\le 0.75$ ($5.98\%$);  
black points satisfy $0.75 <v(q)\le 1$ ($2.77\%$);  
red points satisfy $ v(q)> 1$ ($0.08\%$).
The behaviour of $\G_q$ is the expected one since the red strip 
essentially corresponds with its minimal values, while the minima of $\G^+_q$  seem to be 
less related to $v(q)$; we plan to investigate this phenomenon in the   future.
In Figure \ref{histoEK-EK+}, we also insert two histograms about the 
distribution of the normalised values of  $\G_q$ and $\G^+_q$.
 
In Figures \ref{fig3a}-\ref{fig4} we present the scatter plots 
on $M_q$ and $M_q^\prime :=M_q/\log \log q$.
All the plots were obtained   using GNUPLOT, v.5.2, patchlevel 8,
with the exceptions of the histograms that were obtained using 
Python 3.9 (matplotlib v.3.3.3).

\vfill\eject

\newpage


\vskip 0.5cm
\noindent 
Alessandro Languasco, Luca Righi\\
Universit\`a di Padova,
Dipartimento di Matematica
 ``Tullio Levi-Civita'',\\
Via Trieste 63,
35121 Padova, Italy. \\
{\it e-mail}: alessandro.languasco@unipd.it \\
{\it e-mail}: righi@math.unipd.it 

\end{document}